\documentclass[journal]{new-aiaa}

\usepackage[table]{xcolor}
\usepackage{graphicx}
\usepackage{amsmath,amssymb,bm}
\usepackage{subcaption}
\usepackage{longtable,tabularx}
\usepackage{placeins}
\usepackage{hyperref}
\hypersetup{
colorlinks=true,
linkcolor=brown!90!black,
urlcolor=blue!50,
citecolor=green!50!gray
}

\renewcommand{\arraystretch}{0.8}
\newcommand{\uFTA}{\bm{u}^{\mr{\textsc{fta}}}}

\newcommand{\xinit}{x^{\mr{i}}}
\newcommand{\xfinal}{x^{\mr{f}}}
\newcommand{\BRS}{\mc{B}_{t}(\tau)}
\newcommand{\approxBRS}{\hat{\mc{B}}_{k,t}(\tau)}
\newcommand{\Sigmamk}{\check{\Sigma}^{\mr{m}}_k}
\newcommand{\dotSigmamk}{\dot{\check{\Sigma}}{}^{\mr{m}}_k}
\newcommand{\aactj}{a^{\mr{act}}_j}
\newcommand{\bactj}{b^{\mr{act}}_j}
\newcommand{\cactj}{c^{\mr{act}}_j}
\newcommand{\apskttau}{a^{\mr{ps}}_{k}(t,\tau)}
\newcommand{\bpskttau}{b^{\mr{ps}}_{k}(t,\tau)}
\newcommand{\cpskttau}{c^{\mr{ps}}_{k}(t,\tau)}
\newcommand{\aackjt}{a^{\mr{ac}}_{k,j}(t)}
\newcommand{\backjt}{b^{\mr{ac}}_{k,j}(t)}
\newcommand{\cackjt}{c^{\mr{ac}}_{k,j}(t)}
\newcommand{\aackt}[1]{a^{\mr{ac}}_{k,#1}(t)}
\newcommand{\backt}[1]{b^{\mr{ac}}_{k,#1}(t)}
\newcommand{\cackt}[1]{c^{\mr{ac}}_{k,#1}(t)}

\colorlet{revision}{black}
\newcommand{\rtx}[1]{\textcolor{revision}{#1}}

\newcommand{\behcet}{Beh\c{c}et}

\newcommand{\acikmese}{A\c{c}\i kme\c{s}e}

\newcommand{\zeros}[1]{0_{#1}}
\newcommand{\ones}[1]{1_{#1}}
\newcommand{\eye}[1]{I_{#1}}

\newcommand{\selector}[1]{{}^{#1}\!\!~E}

\newcommand{\tinit}{t^{\mathrm{i}}}
\newcommand{\tfinal}{t^{\mathrm{f}}}
\newcommand{\tsafe}{t^{\mathrm{s}}}
\newcommand{\tburn}{t^{\mr{burn}}}
\newcommand{\tspan}{[\tinit,\tfinal]}
\newcommand{\tsafespan}{[0,\tsafe]}
\renewcommand{\dot}[1]{\overset{\text{\raisebox{-0.5ex}{\large.}}}{#1}}
\newcommand{\derv}[1]{\overset{\raisebox{-1ex}{$\scriptscriptstyle\circ$}}{#1}}
\newcommand{\mr}[1]{\mathrm{#1}}
\newcommand{\mb}[1]{\mathbb{#1}}
\newcommand{\mc}[1]{\mathcal{#1}}
\newcommand{\bR}{\mb{R}}

\newcommand{\sd}[2]{d(#1,#2)}
\newcommand{\projboundary}[2]{{#1}^{\scriptscriptstyle\partial#2}}
\newcommand{\rotmatLVLH}{T^{\scriptscriptstyle\textsc{lvlh}}}
\newcommand{\ctscvx}{{\scalebox{1.1}{\textsc{{\scalebox{0.73}{ct-}}sc{\scalebox{0.73}{vx}}}}}}

\newcommand{\Havoid}{H^{\scriptscriptstyle\mc{A}}}
\newcommand{\havoid}{h^{\scriptscriptstyle\mc{A}}}

\newcommand{\aJtwo}{a^{\scriptscriptstyle{J_2}}}
\newcommand{\aSRP}{a^{\scriptscriptstyle\mr{SRP}}}
\newcommand{\rotmatPA}{T^{\scriptscriptstyle\mr{PA}}}
\newcommand{\Rmoon}{R^{\mr{M}}}
\newcommand{\rsc}{r^{\text{sc}}}
\newcommand{\vsc}{v^{\text{sc}}}
\newcommand{\drsc}{\dot{r}^{\text{sc}}}
\newcommand{\dvsc}{\dot{v}^{\text{sc}}}
\newcommand{\rE}{r^{\text{E}}}
\newcommand{\rS}{r^{ \text{S}}}

\newcommand{\Asc}{A^{\text{sc}}}

\title{{\Large Successive Convexification for Passively-Safe Spacecraft Rendezvous on Near Rectilinear Halo Orbit}}

\author{Purnanand Elango\footnote{Research Scientist; Corresponding author: \texttt{elango@merl.com}},~~Abraham P. Vinod\footnote{Principal Research Scientist}}
\affil{\vspace{0.05cm}Mitsubishi Electric Research Laboratories, Cambridge, MA, 02139, USA}

\author{Kenji Kitamura\footnote{Researcher, Advanced Technology R\&D Center}}
\affil{\vspace{0.05cm}Mitsubishi Electric Corporation, Amagasaki, Hyogo, 661-8661, Japan}

\author{{\behcet}~{\acikmese}\footnote{Professor, Dept. of Electrical Engineering and Computer Science~\&~Dept. of Mechanical Engineering}}
\affil{\vspace{0.05cm}University of California, Berkeley, Berkeley CA, 94720, USA.}

\author{Stefano Di Cairano\footnote{Distinguished Research Scientist},~~Avishai Weiss\footnote{Senior Principal Research Scientist}}
\affil{\vspace{0.05cm}Mitsubishi Electric Research Laboratories, Cambridge, MA, 02139, USA}

\begin{document}
\maketitle
\begin{abstract}
\rtx{We present an optimization-based approach for fuel-efficient spacecraft rendezvous to Gateway, a lunar space station architecture based on a near rectilinear halo orbit (NRHO) around the Earth-Moon L2 point. Our main contribution is a method that ensures safety at all times in the presence of Gaussian uncertainty while incorporating disturbance compensation feedback through correction maneuvers. Specifically, the approach: i) ensures passive safety and satisfies path constraints at all times, ii) satisfies the specifications for critical decision points along the trajectory, and iii) accounts for common uncertainties arising from orbital insertion maneuvers, actuation, and navigation measurements. To account for these uncertainties, we use both chance constraints and disturbance compensation feedback to bound the effects of errors. We leverage sequential convex programming (SCP) together with an isoperimetric reformulation of continuous-time path constraints, including passive safety, to eliminate the risk of inter-sample constraint violations common in existing methods. This enables sparse time discretization that supports efficient offline mission design and analysis and aligns with operational preferences for less frequent control actions. We demonstrate the method on a case study representative of Gateway rendezvous operations.}
\end{abstract}
\clearpage
%
%
\section*{Nomenclature}
{\renewcommand\arraystretch{1.0}
\begin{tabular}{ll}
$\bR,\bR_+$ & Set of real numbers and set of nonnegative real numbers\\
$\bR^n,\bR^{n\times m}$ & Set of $n\times 1$ vector of real numbers and set of real-valued $n\times m$ matrices\\
$\zeros{n},\ones{n}$ & Vector of zeros and vector of ones in $\bR^n$\\
$\eye{n},\zeros{n\times m}$ & Identity matrix in $\bR^{n\times n}$ and matrix of zeros in $\bR^{n\times m}$\\
$(v,w)$ & Concatenation of $v\in\bR^n$ and $w\in\bR^m$ to form a vector in $\bR^{n+m}$\\
$[A~B]$ & Concatenation of $A\in\bR^{n\times m_1}$ and $B\in\bR^{n\times m_2}$ to form a matrix in $\bR^{n\times (m_1+m_2)}$\\ 
$\mr{diag}(v)$ & Diagonal matrix with elements of vector $v$ along the diagonal\\
$\mr{blkdiag}(A,B)$ & Block diagonal matrix formed using square matrices $A$ and $B$\\
$a \le b$ & Elementwise inequalities $a_i\le b_i$, where $a=(a_1,\ldots,a_m)$, $b=(b_1,\ldots,b_m)\in\bR^m$\\
$\|v\|$ & Euclidean norm of vector $v$\\
$\dot{\square},\,\derv{\square}$ & \scalebox{0.8}{$\displaystyle\frac{\mathrel{\raisebox{-3pt}{$\mr{d}\square$}}}{\mathrel{\raisebox{3pt}{$\mr{d}t$}}},\,\frac{\mathrel{\raisebox{-3pt}{$\mr{d}\square$}}}{\mathrel{\raisebox{3pt}{$\mr{d}\tau$}}}$}: time derivatives along the planning and safety horizons, respectively\\[-0.1cm]
$\nabla_m h(x_1,\ldots,x_n)$ & Gradient of $h$ with respect to $m$\textsuperscript{th} argument ($m\le n$), evaluated at $x_1,\ldots,x_n$\\
$\nabla h(x)$ & $\nabla_1 h(x)$\\
$\mc{N}(\mu,\Sigma)$ & Gaussian distribution with mean $\mu$ and covariance $\Sigma$\\
$\bm{x}$ & Boldface symbols denote random variables\\
$\mb{E}(\bm{a})$ & Expected value of random variable $\bm{a}$\\
$\mb{P}(g(\bm{x})\le0)$ & Probability of satisfaction of $g(\bm{x})\le 0$\\
$Q_n$ & Quantile function of the chi-squared distribution with $n$ degrees of freedom
\end{tabular}%
}
%
\section{Introduction}
\rtx{The Gateway~\cite{fuller2024gateway} is a lunar space station architecture based on a near rectilinear halo orbit (NRHO)~\cite{lee2019gateway} around the Earth-Moon L2 point. Visiting and servicing spacecraft will need to rendezvous with the Gateway, which is an important maneuver for many space operations~\cite{flores-abad2014review,ueda2015practical}. For the Gateway rendezvous scenario, a maneuver needs to: i) ensure passive safety at all times, ii) constrain multiple decision points across the trajectory according to specifications~\cite{irsis2019}, iii) remain fuel efficient, and iv) account for uncertainty arising from sources such as orbital insertion maneuvers, actuation, and navigation measurements. Among these requirements, passive safety is especially challenging because it requires the spacecraft's uncontrolled motion (free drift) to remain outside a keep-out zone over a prescribed time interval~\cite{marsillach2020failsafe}---~\!even in the event of an underburn, where the thruster fires for less than the intended duration.}

\rtx{Rendezvous on Keplerian orbits has been extensively studied and successfully executed since the 1960s~\cite{clohessy1960terminal,chamberlin1964gemini}, but it has yet to be fully demonstrated on non-Keplerian orbits such as the NRHO. Two features make NRHO rendezvous particularly challenging. First, because of the orbit's elongated geometry and large spatial scale, rendezvous maneuvers from orbital insertion through proximity operations and docking can span distances exceeding 1000 km~\cite{nakamura2023rendevous}, making fuel-optimal trajectory design numerically challenging. Second, NRHO dynamics deviate from simple two-body motion: the coupled gravitational effects of the Moon and Earth, with smaller contributions from the Sun, produce strongly nonlinear dynamics, particularly near perilune, while higher-order gravitational perturbations, such as the Moon's $J_2$ term, and non-gravitational effects, such as solar radiation pressure (SRP), contribute additional accelerations~\cite{muralidharan2020control}. This contrasts with low Earth orbit (LEO), where the two-body term is typically dominant and lower-order models often suffice. These challenges have motivated significant recent interest in passively-safe rendezvous on NRHO. Prior approaches include center-/unstable-manifold-based design~\cite{bucci2018rendezvous,lizydestrez2019rendezvous}, natural-motion far-range rendezvous in ephemeris models~\cite{blazquez2019rendezvous}, six-degree-of-freedom (6-DoF) close-range rendezvous~\cite{blazquez2021convex}, model predictive control-based guidance~\cite{sanchez2020chance}, reachable-set-based methods~\cite{marsillach2021failsafe}, adjoint methods with ephemeris validation~\cite{bucchioni2021ephemeris}, and linear covariance (LinCov) techniques for fast uncertainty propagation in both gradient-free and gradient-based robust rendezvous trajectory design~\cite{woffinden2022david,goulet2023robust,cavesmith2024angles,cunningham2025robust,margolis2024robust,margolis2024coopt}.}

\rtx{Our goal is to compute a fuel-optimal nominal rendezvous trajectory, from orbital insertion to the start of the proximity phase, while ensuring safety at all times in the presence of uncertainty. Optimization-based methods are natural candidates for this purpose because they can generate trajectories for nonlinear systems subject to complex constraints~\cite{betts2010practical,malyuta2022convex}. They can incorporate safety constraints~\cite{breger2008safe,marsillach2020abortsafe}, temporal-logic specifications~\cite{szmuk2020successive,uzun2024optimization}, and multi-phase dynamics~\cite{olivares2013multiphase,kamath2023real}. These methods have been successfully applied across a wide range of spacecraft operations, including rendezvous~\cite{malyuta2021advances,wang2024survey}, and can be tailored for execution on resource-constrained, radiation-hardened onboard processors~\cite{kamath2023customized,doll2025hardware}. Recent advances have also improved their robustness to numerical ill-conditioning caused by orders-of-magnitude differences in trajectory optimization parameters~\cite{ross2018scaling,chari2024constraint,kamath2025optimal}. Uncertainty and disturbance effects can be incorporated through deterministic reformulations of chance constraints~\cite{blackmore2011chance,caleb2025chance}, covariance steering~\cite{oguri2022stochastic,liu2025optimal,kumagai2025robust}, conditional value-at-risk minimization~\cite{echigo2024dispersion}, and stochastic reachability~\cite{vinod2021abortsafe}. For passive safety, some prior optimization-based approaches explicitly parametrize the free-drift trajectory~\cite{breger2008safe,blazquez2021convex} while others use reachable sets constructed from the keep-out zone~\cite{marsillach2021failsafe}.}

\rtx{Many optimization-based approaches enforce constraints only at a finite set of points, such as discrete time samples in model predictive control or discrete nodes in collocation methods for numerical optimal control~\cite{betts2010practical}. This may lead to inter-sample constraint violations~\cite{dueri2017trajectory} unless the constraints are enforced on a dense grid, which can substantially increase problem size. Related gradient-free approaches, such as particle-swarm and genetic-algorithm methods, iteratively generate and evaluate candidate trajectories. When combined with LinCov, they can account for uncertainty in the cost, constraints, and performance metrics~\cite{woffinden2022david,goulet2023robust,cavesmith2024angles,cunningham2025robust}. However, because gradient-free methods rely on repeated candidate generation and evaluation, they can lead to a large number of trajectory simulations, which can make fine temporal sampling of the cost and constraints computationally expensive. In contrast, gradient-based LinCov approaches embed covariance propagation directly within the optimization and exploit gradient information rather than evolutionary search, which can improve computational efficiency~\cite{margolis2024robust,margolis2024coopt}.}

\rtx{In this work, we address requirements i)--iv) within a gradient-based optimization framework. Specifically, we propose a method for spacecraft rendezvous to the Gateway that ensures passive safety at all times, even in the event of an underburn, and enforces the approach-cone path constraint in continuous time; satisfies the specifications at key decision points along the rendezvous trajectory; and accounts for uncertainty due to NRHO insertion, actuation error, and navigation measurements through chance constraints together with disturbance compensation feedback implemented through correction maneuvers. Our method uses sequential convex programming (SCP) based on the recently proposed continuous-time successive convexification ({\ctscvx}) framework~\cite{elango2025ctscvx}. Within this framework, continuous-time path constraints are converted into integral constraints through an isoperimetric reformulation~\cite{hartl1995survey}, and the passive-safety constraint is modeled using approximate backward reachable sets constructed from the local linearization at each SCP iteration. Because the method is model-based, it is not tied to Gateway or NRHO and can be adapted to rendezvous in low Earth orbit (LEO), geostationary Earth orbit (GEO), highly elliptical Earth orbits, and orbits about Mars.}

\rtx{SCP-based methods have been shown to be flexible and effective for spacecraft trajectory optimization, including rendezvous~\cite{blazquez2021convex,ishizuka2022Asteroid,oguri2022stochastic,malyuta2021advances}, low-thrust transfer design~\cite{oguri2022stochastic}, asteroid rendezvous missions such as Hayabusa2~\cite{ishizuka2022Asteroid}, and close-range 6-DoF rendezvous guidance~\cite{blazquez2021convex}. The present work builds on these developments by focusing on a combination of requirements that, to our knowledge, are not addressed simultaneously in the cited formulations. In particular, prior SCP-based approaches typically address subsets of requirements i)--iv), and path constraints are often enforced at discrete nodes, which may not guarantee satisfaction between nodes. In addition, chance constraints and disturbance compensation feedback for bounding covariance growth are not explicitly included in the cited formulations. Our approach complements this literature by combining continuous-time enforcement of passive-safety and approach-cone constraints with prescribed probabilistic guarantees under multiple sources of uncertainty, while also accounting for the effect of discrete correction maneuvers executed according to an existing, commonly used scheme that modifies thrust only at maneuver points.}

\rtx{The remainder of the paper is organized as follows. Section~\ref{sec:prb-form} formulates the problem, Section~\ref{sec:scp} develops the SCP-based solution approach, Section~\ref{sec:num-results} presents a Gateway case study, and Section~\ref{sec:conclusions} concludes the paper.}
%
\section{Problem Formulation}\label{sec:prb-form}
\begin{figure}[!htpb]
    \centering
    \includegraphics[width=\linewidth]{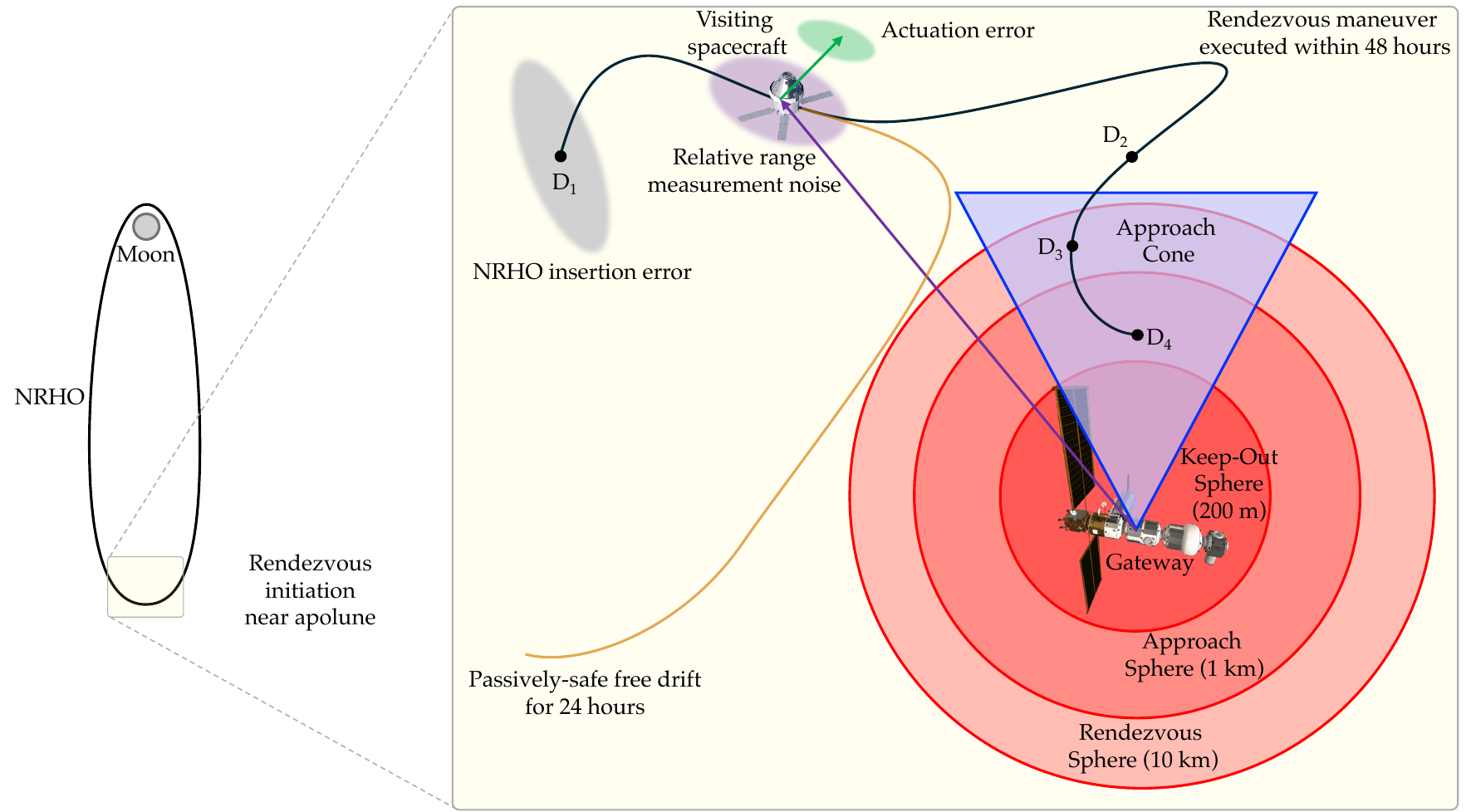}
    \caption{Schematic describing the passively-safe rendezvous maneuver of a visiting spacecraft to the Gateway deployed on an NRHO.}
    \label{fig:nrho-rdv-formulation}
\end{figure}
\rtx{The International Rendezvous System Interoperability Standards (IRSIS)~\cite{irsis2019}, developed by the International Space Station (ISS) member agencies, facilitate collaborative spacecraft operations in cislunar and deep-space environments. We adopt the IRSIS guidelines for spacecraft rendezvous with the Gateway and summarize the maneuver, requirements, and sources of uncertainty in Figure~\ref{fig:nrho-rdv-formulation}. The figure also shows the NRHO region in which rendezvous is initiated, the avoid sets, and the approach cone. The maneuver begins near the NRHO \textit{apolune}, the farthest point on the NRHO from the Moon, after the spacecraft is inserted onto the NRHO, and must reach a hold point 500~m from the Gateway within 48~hr. The hold point lies along the Sun--Gateway axis, and the visiting spacecraft must approach the Gateway within a cone about this axis so that the Gateway is illuminated and visible to the spacecraft's optical navigation sensors. The maneuver is divided into three phases with four decision points: D\textsubscript{1}, where rendezvous begins after NRHO insertion and 24-hr passive safety must be maintained with respect to the rendezvous sphere; D\textsubscript{2}, where entry into the rendezvous sphere begins and 24-hr passive safety must be maintained with respect to the approach sphere; D\textsubscript{3}, where entry into the approach sphere begins and 24-hr passive safety must be maintained with respect to the keep-out sphere; and D\textsubscript{4}, where entry into the keep-out sphere begins for the proximity phase. At D\textsubscript{1}, the relative distance between the spacecraft and the Gateway can range from 400 to 2000~km~\cite{nakamura2023rendevous}; at D\textsubscript{2}, D\textsubscript{3}, and D\textsubscript{4}, the spacecraft is required to lie outside the rendezvous sphere, the approach sphere, and the keep-out sphere, respectively. We consider three sources of stochastic uncertainty in the guidance design: i) NRHO insertion error, ii) spacecraft thruster actuation error, and iii) noise in relative range and range-rate measurements. We refer the reader to~\cite[Chap.~4]{fehse2003automated} for a detailed discussion of safety considerations and uncertainties in Keplerian-orbit rendezvous. The same underlying elements---hold points, avoid sets, passive safety, and the approach cone---also apply to NRHO rendezvous.}

Unless otherwise noted, our modeling assumptions and numerical specifications (avoid sets, decision-point windows, approach-cone/range parameters, and uncertainty statistics) are based on IRSIS~\cite{irsis2019} and prior numerical studies of rendezvous to the Gateway~\cite{nakamura2023rendevous}.

A fuel-efficient rendezvous maneuver that \rtx{satisfies these requirements} while accounting for uncertainty can be obtained by solving a stochastic optimal control problem. The construction of this problem is detailed in the following subsections.
%
\subsection{Spacecraft Relative Motion Model}
Consider a visiting spacecraft in the vicinity of the Gateway. The (uncontrolled) natural motion of the Gateway in the Moon-centered J2000 inertial frame~\cite[Sec. 3.7]{vallado2013fundamentals} is described by
\begin{align}
    \dot{z}(t) ={} & \check{f}(t,z(t)),\label{cislunar-sc-dyn}
\end{align}
where $t\in\bR_+$ and $z(t)\in\bR^6$. The natural motion, also referred to as \textit{free drift}, is governed by gravitational acceleration from the Moon, Earth, and Sun, along with higher-order effects such as solar radiation pressure (SRP) and the Moon's $J_2$ gravitational coefficient. The full dynamical model is given in Appendix~\ref{app:sc-eom}. The state trajectory $x$ of the visiting spacecraft, expressed relative to the Gateway~\cite{lee2019gateway}, satisfies
\begin{align}
    \dot{x}(t) ={} & f(t,x(t),u(t)) \triangleq \check{f}(t,x(t)+z(t)) - \check{f}(t,z(t)) + \begin{bmatrix}\zeros{3\times 3}\\[0.05cm]\eye{3}\end{bmatrix}u(t),\label{visit-sc-dyn}
\end{align}
where $u$ is a piecewise-continuous control input acting on the spacecraft. \rtx{It represents either an acceleration or a sequence of velocity impulses; the latter is represented in continuous time using Dirac delta functions.} Spacecraft are controlled by chemical or electric thrusters that are activated for a \rtx{finite duration determined by} the required momentum change. These \rtx{finite-duration control actions can} be modeled using zero-order hold (ZOH)~\cite[Chap. 3]{wittenmark2002computer} or as fixed-duration pulses, which we refer to as \textit{finite-burn pulses} (FBPs). When the burn duration is small relative to the time scale of the trajectory, an impulsive model is often sufficiently accurate. We adopt the mass-normalized formulation in~\eqref{visit-sc-dyn}, which is convenient for design. \rtx{The state and control input dimensions are} $n^x = 6$ and $n^u = 3$, respectively.
\subsection{Path Constraints}
The visiting spacecraft must satisfy passive-safety and approach-cone path constraints along its rendezvous trajectory to the Gateway. The passive-safety constraint requires that the \rtx{spacecraft's free-drift trajectory}, starting at time $t$, not enter an open set $\mc{A}$, referred to as an \textit{avoid set}, over a time horizon of length $\tsafe$. \rtx{This requirement can be stated precisely} using backward reachable sets (BRS). Given $t,\tau\in\bR_+$, the BRS of $\mc{A}$ over the horizon $[t,t+\tau]$ is the set of all states at time $t$ that can reach $\mc{A}$ at time $t+\tau$ under natural motion (i.e., free-drift dynamics)
\begin{align}
    \BRS \triangleq{} \left\{ x(t)\in\bR^{n^x}\,\middle|~\begin{array}{l}  x(t+\tau)\in\mc{A} \\\displaystyle\frac{\mr{d}x(\gamma)}{\mr{d}\gamma} = f(\gamma,x(\gamma),\zeros{n^u}),~\gamma\in[t,t+\tau]\end{array}\right\}.
\end{align}%
We refer to~\cite[Sec. 14.2.1.3]{lavalle2006planning} and~\cite[Def. 10.5]{borrelli2017predictive} for related notions of BRS. Given a state trajectory $x$ satisfying~\eqref{visit-sc-dyn}, the passive-safety constraint at time $t$ can be stated as
\begin{align}
    x(t) \notin \bigcup_{\tau\in\tsafespan}\!\BRS,\label{ps-cnstr-BRS-naive}
\end{align}
where $\tsafespan$ is referred to as the \textit{safety horizon}. \rtx{Thus, requiring the spacecraft's free-drift trajectory, starting at time $t$, to remain outside $\mc{A}$ at all times over $[t,t+\tsafe]$ is equivalent to requiring the state $x(t)$ to lie outside the collection of BRS in~\eqref{ps-cnstr-BRS-naive}.}

\rtx{Equivalently,~\eqref{ps-cnstr-BRS-naive} can be written} in terms of the signed distance from a point to a set~\cite[Ex. 8.5(b)]{boyd2004convex}
\begin{align}
    \rtx{\sd{x(t)}{\BRS}} > 0, \quad \forall\,\tau\in\tsafespan.\label{ps-cnstr-sd-naive}  
\end{align}
\rtx{Strict inequalities are generally difficult to impose in numerical optimization, so we relax} the strict inequality in~\eqref{ps-cnstr-sd-naive} to a non-strict inequality, assuming that the BRS are open sets (a point $x$ on the boundary of an open set $\mc{B}$ satisfies $\rtx{\sd{x}{\mc{B}}} = 0$). In the proposed solution method, the BRS are iteratively approximated with polyhedral sets by linearizing~\eqref{visit-sc-dyn}. \rtx{Because $\mc{A}$ is open, the approximate BRS are also open sets (see Section~\ref{subsec:approx-BRS} for details). We use the signed distance function because it is differentiable almost everywhere for sets with piecewise smooth boundaries and can be computed efficiently for polyhedral sets (see Appendix~\ref{app:signed-distance}).}

The approach-cone constraint requires the spacecraft to approach the origin, i.e., the Gateway, from within a cone with specified axis $e^{\mr{ac}}$ and half-angle $\theta^{\mr{ac}}$. Given a state trajectory $x$ satisfying~\eqref{visit-sc-dyn}, \rtx{this requirement} can be stated as the second-order cone constraint
\begin{align}
    \cos\theta^{\mr{ac}}\|\selector{r}x(t)\|\le (e^{\mr{ac}})^\top\selector{r}x(t),\label{ac-soc}
\end{align}
where $\selector{r}$ selects the position coordinates from a vector in $\bR^{n^x}$. \rtx{The proposed method requires gradient information, whereas the Euclidean norm in~\eqref{ac-soc} is not differentiable at the origin. We therefore reformulate~\eqref{ac-soc} in an equivalent, continuously differentiable form by squaring both sides and additionally requiring the right-hand side of~\eqref{ac-soc} to be nonnegative. To this end, define $g^{\mr{ac}}:\bR^{n^x}\to\bR^2$ as}
\begin{align}
    g^{\mr{ac}}(x(t)) \triangleq 
    \begin{bmatrix} g^{\mr{ac}}_1(x(t)) \\[0.1cm] g^{\mr{ac}}_2(x(t)) \end{bmatrix}
    \triangleq 
    \begin{bmatrix}
        \|\cos\theta^{\mr{ac}}\,\selector{r}x(t)\|^2 - \big((e^{\mr{ac}})^\top\selector{r}x(t)\big)^2\\
        -(e^{\mr{ac}})^\top\selector{r}x(t)
    \end{bmatrix},
\end{align}
and express the approach-cone constraint as
\begin{align}
    g^{\mr{ac}}(x(t)) \le \zeros{2}. \label{ac-diff}
\end{align}%
%
%
\subsection{Stochastic Optimal Control Problem}
\begin{figure}[!htp]
    \centering
    \includegraphics[width=\linewidth]{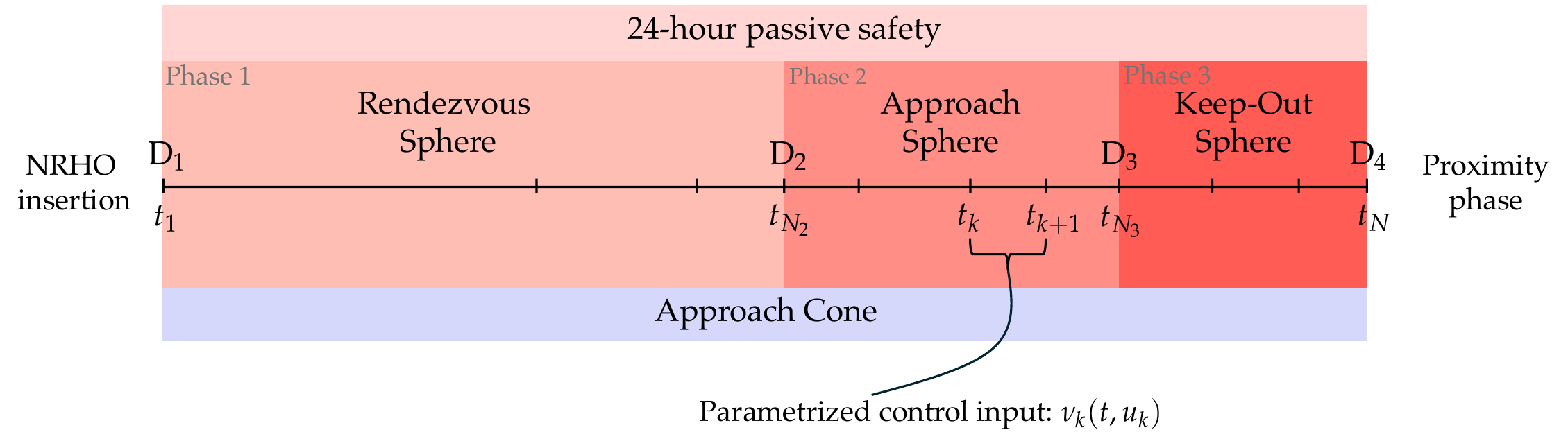}
    \caption{Schematic timeline of the rendezvous maneuver showing the decision points and the phases over which the path constraints are imposed.}
    \label{fig:maneuver-timeline}
\end{figure}

We now formulate a stochastic optimal control problem for the rendezvous maneuver illustrated in Figure~\ref{fig:nrho-rdv-formulation}. The maneuver is executed over the time interval $\tspan$, which we refer to as the \textit{planning horizon}. We assume that actuation and navigation measurements \rtx{occur on a grid of $N$ time nodes in $\tspan$}, denoted by $\tinit=t_1,\ldots,t_N=\tfinal$. The decision points D\textsubscript{1}, D\textsubscript{2}, D\textsubscript{3}, and D\textsubscript{4} in Figure~\ref{fig:nrho-rdv-formulation} occur at times $t_1$, $t_{N_2}$, $t_{N_3}$, and $t_N$, respectively, where $1 < N_2 < N_3 < N$. \rtx{We denote the avoid set for phase $j$ by $\mc{A}_j$, for $j=1,2,3$; accordingly, $\BRS$ denotes the BRS of $\mc{A}_1$ for $t\in[\tinit,t_{N_2})$, of $\mc{A}_2$ for $t\in[t_{N_2},t_{N_3})$, and of $\mc{A}_3$ for $t\in[t_{N_3},\tfinal]$. For each $k=1,\ldots,N-1$, the control input signal over $[t_k,t_{k+1}]$ is parameterized as an impulse, a ZOH input, or an FBP input. Specifically, for $t\in[t_k,t_{k+1}]$, the control input is given by} the parametrization function $\nu_k : [t_k,t_{k+1}] \times\bR^{n^u} \to \bR^{n^u}$ as
\begin{align}
    \nu_k(t,u_k) \triangleq \begin{cases}\delta(t-t_k)u_k &\text{Impulse},\\
    u_k &\text{ZOH},\\
    \left.\begin{array}{ll}
        u_k &~~\text{if } t \le t_k+\tburn\\
        \zeros{n^u} &~~\text{otherwise}\end{array}
    \right\} & \text{FBP}, 
    \end{cases}    
\end{align}
for $u_k\in\bR^{n^u}$, where $\delta$ is the Dirac delta function and $\tburn$ is the \rtx{thruster burn duration}. We treat $u_1,\ldots,u_{N-1}$ as decision variables in the stochastic optimal control problem and refer to them as \textit{nominal} control inputs. For the remainder of the development (i.e., the solution approach in Section~\ref{sec:scp}), the time grid $t_1,\ldots,t_N$ \rtx{is treated as given}. Operational considerations for choosing this grid are discussed in Section~\ref{sec:num-results} in the context of a realistic numerical case study. \rtx{Figure~\ref{fig:maneuver-timeline} shows the schematic timeline of the rendezvous maneuver and the phases over which the path constraints are imposed. The differently shaded red regions indicate the phase-dependent avoid sets for which 24-hr passive safety is enforced, and the blue region indicates the approach-cone constraint imposed throughout the maneuver.}

\rtx{For each $k=1,\ldots,N-1$, the flow map obtained by integrating~\eqref{visit-sc-dyn} with the parametrized control input signal is denoted by $F_k:\bR^{n^x}\times\bR^{n^u} \to \bR^{n^x}$ and defined as}
\begin{align}
    F_k(x_k,u_k) \triangleq{} & x_k + \int_{t_k}^{t_{k+1}}f(t,\check{x}_k(t),\nu_k(t,u_k))\mr{d}t,    
\end{align}
for $x_k\in\bR^{n^x}$, where \rtx{the} state trajectory $\check{x}_k$ satisfies~\eqref{visit-sc-dyn} over $[t_k,t_{k+1}]$ with initial condition $x_k$ and control input signal $t\mapsto \nu_k(t,u_k)$.

The stochastic optimal control problem for the three-phase minimum-fuel rendezvous maneuver is then
\begin{subequations}
\begin{align}
\underset{u_k}{\mr{minimize}}~&~\sum_{k=1}^{N-1} \alpha_k\|u_k\|, & & \label{stoc-ocp:obj}\\
\mr{subject~to}~&~\bm{x}_{k+1} = F_k(\bm{x}_k,u_k+\bm{\mu}_k), & & \hspace{0cm}\hphantom{\tau\in\tsafespan,~}\,t\in[t_k,t_{k+1}],~k=1,\ldots,N-1,\label{stoc-ocp:disc-dyn}\\
&~\mb{P}(\rtx{\sd{\check{\bm{x}}_k(t)}{\BRS}} \ge 0)\ge \beta^{\mr{ps}}, & & \hspace{0cm}\tau\in\tsafespan,~t\in[t_{k},t_{k+1}],~k=1,\ldots,N-1,\label{stoc-ocp:ps-cnstr}\\
&~\mb{P}(g^{\mr{ac}}(\check{\bm{x}}_k(t)) \le \zeros{2})\ge \beta^{\mr{ac}}, & & \hspace{0cm}\hphantom{\tau\in\tsafespan,~}\,t\in[t_{k},t_{k+1}],~k=1,\ldots,N-1,\label{stoc-ocp:ac-cnstr}\\
&~\mb{P}(u_k + \bm{\mu}_k \in \mc{U})\ge \beta^{\mr{act}},~\bm{\mu}_k\sim\mc{N}(\zeros{n^u},\Sigma^{\mr{act}}), & &\hspace{0cm}\hphantom{\tau\in\tsafespan,~t\in[t_{k},t_{k+1}],~}\,k=1,\ldots,N-1, \label{stoc-ocp:ctrl-cnstr}\\
&~\|\selector{r}\,\mb{E}(\bm{x}_{N_{j}})\| \le a^+_{j},~(e^{\mr{ac}})^\top\selector{r}\,\mb{E}(\bm{x}_{N_{j}}) \ge a^-_{j}, & &\hspace{0cm}\hphantom{\tau\in\tsafespan,~t\in[t_{k},t_{k+1}],~}\,~j=2,3,\label{stoc-ocp:phase-transition}\\
&~\bm{x}_1 \sim \mc{N}(\xinit,\Sigma^{\mr{i}}),~\mb{E}\big(\bm{x}_N\big) = \xfinal.\label{stoc-ocp:bc}
\end{align}\label{stoc-ocp}%
\end{subequations}
The initial state $\bm{x}_1$ at D\textsubscript{1} is a random variable with distribution $\mc{N}(\xinit,\Sigma^{\mr{i}})$, and the nominal final state at D\textsubscript{4} is denoted by $\xfinal$. For each $k=1,\ldots,N-1$, the spacecraft trajectory $\check{\bm{x}}_k$ satisfies~\eqref{visit-sc-dyn} over $[t_k,t_{k+1}]$ with initial condition $\bm{x}_k$ and control input signal $t\mapsto\nu_k(t,u_k+\bm{\mu}_k)$. Here, the nominal control input $u_k$ is perturbed by the random variable $\bm{\mu}_k\sim\mc{N}(\zeros{n^u},\Sigma^{\mr{act}})$. The coefficient $\alpha_k$ in~\eqref{stoc-ocp:obj} depends on the control parametrization: for ZOH, $\alpha_k = t_{k+1} - t_k$, whereas for impulse and FBP, $\alpha_k = 1$. The chance constraints~\eqref{stoc-ocp:ps-cnstr},~\eqref{stoc-ocp:ac-cnstr}, and~\eqref{stoc-ocp:ctrl-cnstr} are enforced with confidence levels $\beta^{\mr{ps}}$, $\beta^{\mr{ac}}$, and $\beta^{\mr{act}}$, respectively. The set of admissible control inputs $\mc{U}$ is \rtx{assumed to be compact and convex. The decision-point requirements are enforced through~\eqref{stoc-ocp:phase-transition}--\eqref{stoc-ocp:bc}; in particular,~\eqref{stoc-ocp:phase-transition}} ensures that the expected spacecraft positions at D\textsubscript{2} and D\textsubscript{3} lie in the intersection of a ball and a halfspace.

Satisfying the constraints in~\eqref{stoc-ocp} can be challenging \rtx{without} a feedback mechanism to mitigate the effect of uncertainty. To address this, the proposed solution approach incorporates disturbance compensation feedback through correction maneuvers based on state measurements. The path chance constraints~\eqref{stoc-ocp:ps-cnstr} and~\eqref{stoc-ocp:ac-cnstr} must be satisfied in continuous time, but a discrete-time model of the spacecraft dynamics is sufficient in~\eqref{stoc-ocp:disc-dyn} because the insertion, actuation, and measurement uncertainties manifest only at discrete time instants in the planning horizon. The discrete-time dynamics model in~\eqref{stoc-ocp:disc-dyn} is exact, \rtx{which avoids the need for a stochastic differential-equation model and keeps the subsequent development straightforward.}
%
\section{Sequential Convex Programming-Based Solution Approach}\label{sec:scp}
The proposed SCP-based approach for solving~\eqref{stoc-ocp} \rtx{consists of five key steps: i) reformulating the chance constraints, ii) incorporating disturbance compensation feedback through correction maneuvers} based on relative range measurements, iii) approximating BRS for the nonlinear spacecraft dynamics, iv) reformulating path constraints in isoperimetric form, and v) iteratively solving and refining a convexified and relaxed optimal control problem.

In the remainder of the discussion, reference states and control inputs are \rtx{denoted by an overbar (e.g., $\bar{x}$ and $\bar{u}$). For each $k=1,\ldots,N-1$, a discrete-time quantity $\square_k$ at time $t_k$ may be associated with a continuous-time signal over $[t_k,t_{k+1}]$, denoted by $\check{\square}_k(t)$ for $t\in[t_k,t_{k+1}]$. For example, $\check{\bar{x}}_k(t)$ denotes the continuous-time reference state trajectory over $[t_k,t_{k+1}]$ associated with the discrete-time reference state $\bar{x}_k$ at $t_k$.}
%
\subsection{State Covariance Evolution}\label{subsec:cov-evolve-feedback}
In the proposed approach, \rtx{the problem is linearized and convexified about} reference states and control inputs that are updated at each iteration. Consider a sequence of reference states $\bar{x}_1,\ldots,\bar{x}_N$ and control inputs $\bar{u}_1,\ldots,\bar{u}_{N-1}$. Let $k=1,\ldots,N-1$. 
%
%
For any $x_k\in\bR^{n^x}$ and $u_k\in\bR^{n^u}$, the first-order expansion of $F_k$ is denoted by
\begin{align}
    F_k(x_k,u_k) \approx A_kx_k + B_ku_k + c_k, \label{eq:lin-Fk}   
\end{align}%
where $A_k$, $B_k$, and $c_k$ are obtained by solving
\begin{subequations}
\begin{align}
    \dot{\check{\bar{x}}}_k(t) ={} & f(t,\check{\bar{x}}_k(t),\nu_k(t,\bar{u}_k)), & & t\in[t_k,t_{k+1}],\\
    \dot{\Phi}(t,t_k) ={} & \check{A}_k(t) \Phi(t,t_k), & & t\in[t_k,t_{k+1}],\\
    \dot{\Psi}(t,t_k) ={} & \check{A}_k(t) \Psi(t,t_k) + \check{B}_k(t)\nabla_2 \nu_k(t,\bar{u}_k), & & t\in[t_k,t_{k+1}],\\
    \dot{\Theta}(t,t_k) ={} & \check{A}_k(t) \Theta(t,t_k) + \check{c}_k(t), & & t\in[t_k,t_{k+1}],\\
    \check{\bar{x}}_k(t_k) ={}& \bar{x}_k,\\
    \Phi(t_k,t_k) ={}& \eye{n^x},\\
    \Psi(t_k,t_k) ={}& \zeros{n^x\times n^u},\\
    \Theta(t_k,t_k) ={}& \zeros{n^x},
\end{align}\label{eq:fo-sensitivity}%
\end{subequations}
with 
\begin{align*}
\check{A}_k(t) \triangleq{} & \nabla_2 f(t,\check{\bar{x}}_k(t),\nu_k(t,\bar{u}_k)),\\
\check{B}_k(t) \triangleq{} & \nabla_3 f(t,\check{\bar{x}}_k(t),\nu_k(t,\bar{u}_k)),\\
\check{c}_k(t) \triangleq{} & f(t,\check{\bar{x}}_k(t),\nu_k(t,\bar{u}_k)) - \check{A}_k(t)\check{\bar{x}}_k(t) - \check{B}_k(t)\nu_k(t,\bar{u}_k).
\end{align*}
Then, $A_k = \Phi(t_{k+1},t_k)$, $B_k = \Psi(t_{k+1},t_{k})$, and $c_k = \Theta(t_{k+1},t_k)$. \rtx{For the FBP parametrization, $\nabla_2\nu_k(t_k+\tburn,\bar{u}_k)$ is not defined at the switching time. Therefore, the integration of~\eqref{eq:fo-sensitivity} over $[t_k,t_{k+1}]$ is split into $[t_k,t_k+\tburn]$ and $[t_k+\tburn,t_{k+1}]$.} In the latter interval, we set $\nabla_2\nu_k(t_k+\tburn,\bar{u}_k) = \zeros{n^u}$.

Given nominal control inputs $u_1,\ldots,u_{N-1}$, let the corresponding nominal states $x_1,\ldots,x_N$ satisfy the linearized discrete-time spacecraft dynamics~\eqref{eq:lin-Fk},
\begin{align}
    x_{k+1} ={} A_kx_k + B_ku_k + c_k,\label{eq:lin-disc-dyn}
\end{align}
for $k=1,\ldots,N-1$, with $x_1 = \xinit$. Linearizing~\eqref{stoc-ocp:disc-dyn} about this nominal trajectory yields
\begin{align}
\bm{x}_{k+1} = A_k\bm{x}_k + B_k(u_k + \bm{\mu}_k) + c_k,\label{eq:lin-disc-dyn-rand}
\end{align}
with $\bm{x}_1 \sim \mc{N}(\xinit,\Sigma^{\mr{i}})$. Under the Gaussian approximation adopted here, $\bm{x}_k\sim\mc{N}(x_k,\Sigma_k)$ for each $k=1,\ldots,N-1$, where the covariance satisfies
\begin{align}
\Sigma_{k+1} ={} & A_k\Sigma_kA_k^\top + B_k\Sigma^{\mr{act}}B_k^\top,\label{eq:disc-cov-dyn}
\end{align}
with initial condition $\Sigma_1 = \Sigma^{\mr{i}}$. {\color{revision} In open loop, that is, under the nominal control input alone, an unstable or marginally stable system such as~\eqref{visit-sc-dyn} can exhibit substantial covariance growth over time, for example when $A_k$ in~\eqref{eq:disc-cov-dyn} has eigenvalues greater than one. As a result, satisfying constraints reliably under uncertainty becomes difficult. In long-horizon spacecraft operations, covariance growth is often controlled by applying correction maneuvers at discrete decision times to steer the spacecraft back toward the guidance trajectory in the presence of disturbances. These correction maneuvers effectively reduce the state dispersion and keep the covariance at levels compatible with the imposed constraints. A commonly used approach to compute such corrections is the fixed-time-of-arrival (FTA) guidance law, which has substantial heritage in space missions since the Apollo era~\cite{mclean1962optimal}. Here we incorporate the effect of disturbance compensation feedback through FTA correction maneuvers, applied at the same discrete maneuver times as the guidance inputs, into the guidance problem. This requires three steps: 1) represent FTA correction maneuvers in state-feedback form; 2) account for the resulting closed-loop dynamics in the covariance propagation; and 3) quantify the actual worst-case covariance growth under FTA corrections.

Let the measured state at time $t_k$ be $\bm{x}^{\mr{m}}_k \triangleq \bm{x}_k+\bm{\zeta}_k$, where $\bm{\zeta}_k\sim\mc{N}(\zeros{n^x},\Sigma_k^{\mr{rr}})$ is the measurement noise. In the FTA approach, the correction maneuver at time $t_k$ is chosen so that, in the absence of measurement and actuation noise, the position error at time $t_{k+1}$ is driven to zero. Accordingly, the FTA correction maneuver $\uFTA_k$ can be represented in linear state-feedback form with gain $K_k$ satisfying
\begin{align} 
		\selector{r}(A_k(\bm{x}^{\mr{m}}_k-x_k)+B_k\uFTA_k) = 
		\zeros{3}\iff  	\uFTA_k= -(\selector{r}B_k)^{-1}\selector{r}A_k(\bm{x}^{\mr{m}}_k-x_k) 
		\implies 
		K_k = 
		-(\selector{r}B_k)^{-1}\selector{r}A_k.
		\label{gain-compute}
	\end{align}
The inverse in~\eqref{gain-compute} requires $\selector{r}B_k$ to have full rank, i.e., independent first-order position authority in all three control directions. This condition holds for the three-axis actuation and maneuver intervals used here. From~\eqref{eq:lin-disc-dyn-rand}, the measured-state evolution, for $k=1,\ldots,N-1$, is
\begin{align}
		\bm{x}^{\mr{m}}_{k+1} = A_k\bm{x}^{\mr{m}}_k + B_k(u_k + \uFTA_k + \bm{\mu}_k) + 
		\bm{\zeta}_{k+1} - A_k\bm{\zeta}_k + c_k, \label{eq:lin-disc-dyn-rand-meas}        
\end{align}
Using the feedback representation~\eqref{gain-compute}, the applied control input is $u_k + K_k(\bm{x}^{\mr{m}}_k-x_k)$, i.e., the nominal guidance input plus the FTA correction maneuver. Combining~\eqref{eq:lin-disc-dyn},~\eqref{eq:lin-disc-dyn-rand-meas}, and~\eqref{gain-compute} gives
\begin{align}
    \bm{x}^{\mr{m}}_{k+1} - x_{k+1} = A^{\mr{cl}}_k(\bm{x}^{\mr{m}}_k - x_k) + Y_k\begin{bmatrix} \bm{\zeta}_{k+1} \\ \bm{\zeta}_k \\ \bm{\mu}_k \end{bmatrix}\!, 
\end{align}
where $A^{\mr{cl}}_k \triangleq A_k + B_kK_k$ and $Y_k\triangleq 
[\eye{n^x}~\!-\!A_k~B_k]$. 
We emphasize that the disturbance compensation feedback is implemented through FTA correction maneuvers applied only at the same discrete instants as the guidance maneuvers, rather than through continuous-time feedback.}

\begin{figure}[!htp]
    \centering
    \includegraphics[width=0.7\linewidth]{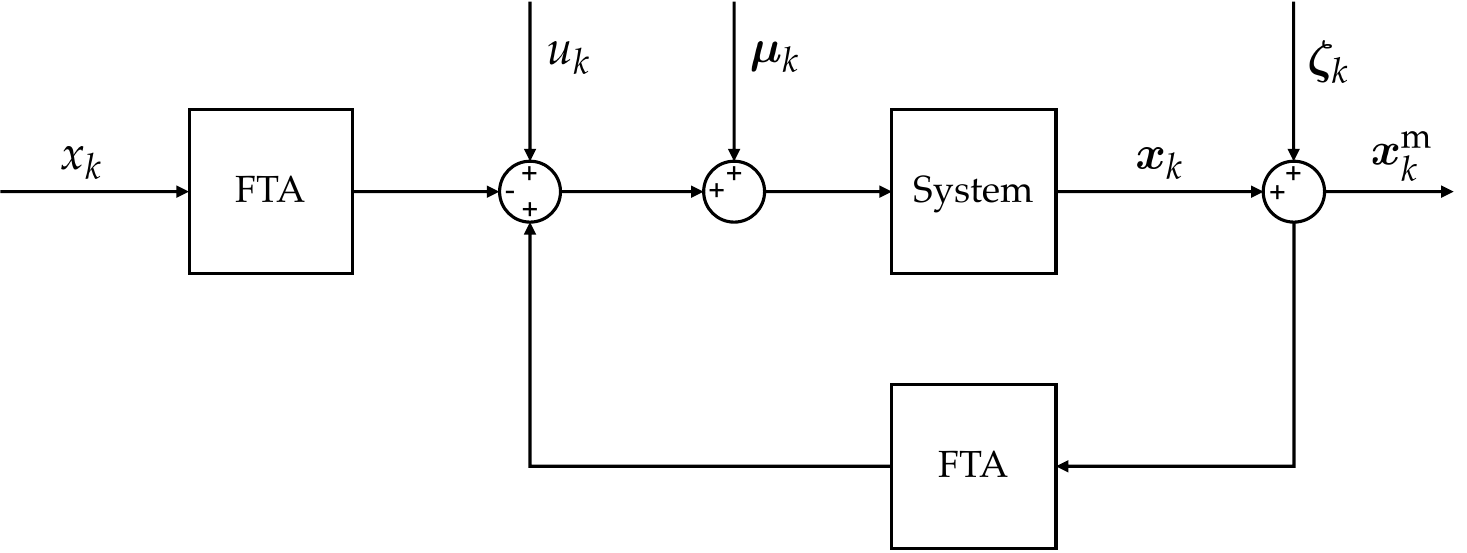}
    \caption{Block diagram of system~\eqref{eq:lin-disc-dyn-rand-meas} with actuation error, measurement noise, and \rtx{disturbance compensation feedback}.}
    \label{fig:blk-cl-sys}
\end{figure}

Since $\bm{\zeta}_1,\ldots,\bm{\zeta}_N$, $\bm{\mu}_1,\ldots,\bm{\mu}_{N-1}$, and $\bm{x}_1$ are independent random variables, $\bm{x}_1^{\mr{m}}\sim\mc{N}(\xinit,\Sigma_1^{\mr{m}})$ with $\Sigma^{\mr{m}}_1 = \Sigma_1 + \Sigma_1^{\mr{rr}}$. Then, for each $k=1,\ldots,N-1$, the covariance of the measured state satisfies
\begin{align}
    \Sigma^{\mr{m}}_{k+1} ={} & A_k^{\mr{cl}}\Sigma^{\mr{m}}_{k}A_k^{\mr{cl}}{}^\top + Y_k\Omega_kY_k^\top,\label{eq:disc-cov-dyn-meas}
\end{align}
where $\bm{x}^{\mr{m}}_k\sim\mc{N}(x_k,\Sigma^{\mr{m}}_k)$ and $\Omega_k \triangleq 
\mr{blkdiag}(\Sigma^{\mr{rr}}_{k+1},\Sigma^{\mr{rr}}_{k},\Sigma^{\mr{act}})$. \rtx{Although we focus here on the case in which FTA is used, and hence obtain $K_k$ as a feedback representation of FTA, the proposed approach remains compatible with any linear feedback law. Figure~\ref{fig:blk-cl-sys} shows the block diagram of the closed-loop system with the feedback structure described thus far.}

For a nonlinear dynamical system such as~\eqref{visit-sc-dyn}, the state distribution over time is generally non-Gaussian, even if the initial state and all injected uncertainties are Gaussian. However, numerical results in~\cite{gupta2025modified} \rtx{show that, for most locations along the NRHO, the spacecraft state distribution is well approximated by a Gaussian when both the initial state and injected uncertainties are Gaussian.} For simplicity, we model actuation error $\bm{\mu}_k$ as a zero-mean Gaussian dispersion chosen conservatively to subsume plausible biases. More realistic, bias-aware models (e.g., the Gates model~\cite{gates1963simplified,geller2006linear}) can be incorporated without altering the proposed solution approach. \rtx{A constant or slowly varying actuation or navigation bias simply requires adjusting the propagated mean~\eqref{eq:lin-disc-dyn-rand-meas} and covariance~\eqref{eq:disc-cov-dyn-meas}, leaving the continuous-time chance-constraint formulation unchanged.} Likewise, model mismatch can be captured as Gaussian process noise with covariance $Q_k$ in~\eqref{eq:disc-cov-dyn-meas}.
\subsection{Reformulation of Chance Constraints}\label{subsec:chance-cnstr-reform}
\rtx{Next, we provide deterministic reformulations of the control input chance constraint~\eqref{stoc-ocp:ctrl-cnstr} and the linearized passive-safety~\eqref{stoc-ocp:ps-cnstr} and approach-cone~\eqref{stoc-ocp:ac-cnstr} path chance constraints.} We adopt a standard reformulation based on ellipsoidal confidence sets~\cite[Rem. 1]{lew2020chance}, which requires the feasible set of each constraint to be polyhedral.

We begin with the control input constraints. If the admissible set $\mc{U}$ is an arbitrary compact convex set, a polytopic inner approximation can be computed efficiently~\cite{ben-tal1999polyhedral,eren2015constrained}. If $\mc{U}$ is a second-order cone or a Euclidean ball, the method in~\cite[Sec.~IV.A.3]{kumagai2025robust} applies. Both approaches introduce inherent conservatism. In this work, we adopt the former approach and assume that the admissible control inputs lie within a polytope, i.e., $\mc{U} = \{u\in\bR^{n^u}\,|\,(\aactj)^\top u \le \bactj,~j=1,\ldots,n^{\mr{act}}\}$. For each $k=1,\ldots,N-1$, the closed-loop control input signal \rtx{over} $[t_k,t_{k+1}]$ is determined by the random variable $\bm{u}^{\mr{cl}}_k = u_k+\bm{\mu}_k + K_k(\bm{x}^{\mr{m}}_k-x_k)$, with distribution $\mc{N}(u_k,\Sigma^{\mr{act}}+K_k\Sigma^{\mr{m}}_kK_k^\top)$. The chance constraint~\eqref{stoc-ocp:ctrl-cnstr} imposed on $\bm{u}_k^{\mr{cl}}$ can be reformulated conservatively using ellipsoidal confidence sets as
\begin{align}
    (\aactj)^\top u_k \le \bactj - \cactj,~j=1,\ldots,n^{\mr{act}},\label{ctrl-chance-reform}
\end{align}
where 
\begin{align}
    \cactj \triangleq  \sqrt{Q_{{n^u}}(\beta^{\mr{act}})(\aactj)^\top(\Sigma^{\mr{act}} + K_k\Sigma^{\mr{m}}_kK_k^\top)\aactj},
\end{align}%
and $Q_{{n^u}}$ is the quantile function of the chi-squared distribution with $n^u$ degrees of freedom~\cite[Cor. 2]{hewing2018on}. Satisfaction of~\eqref{ctrl-chance-reform} implies that
\begin{align*}
    \mb{P}\big((\aactj)^\top \bm{u}^{\mr{cl}}_k \le \bactj,~j=1,\ldots,n^{\mr{act}}\big) \ge \beta^{\mr{act}}.
\end{align*}

To reformulate the passive-safety and approach-cone path chance constraints, we use continuous-time trajectories of the state mean and covariance that align with the discrete-time evolution in~\eqref{eq:lin-disc-dyn} and~\eqref{eq:disc-cov-dyn-meas}. For each $k=1,\ldots,N-1$, the solution to
\begin{subequations}
\begin{align}
    \dot{\check{x}}{}_k(t) ={} & \check{A}_k(t)\check{x}_k(t) + \check{B}_k(t)\nu_k(t,u_k) + \check{c}_k(t), & & t\in[t_k,t_{k+1}],\label{eq:lin-ct-dyn}\\
    \dotSigmamk(t) ={} & \check{A}_k(t)\Sigmamk(t) + \Sigmamk(t)\check{A}_k(t){}^\top, & & t\in[t_{k},t_{k+1}],\label{eq:ol-cov-evolve}\\
    \check{x}_k(t_k) ={} & x_k, \\
    \Sigmamk(t_k) ={} & (\eye{n^x} + A_k^{-1}B_k K_k)\Sigma^{\mr{m}}_k(\eye{n^x} + A_k^{-1}B_k K_k)^\top + A_k^{-1}Y_k \Omega_kY_k^\top A_k^{-\!\top},\label{eq:ol-cov-init}
\end{align}%
\end{subequations}
provides $\check{x}_k(t_{k+1}) = x_{k+1}$ \rtx{and} $\Sigmamk(t_{k+1}) = \Sigma^{\mr{m}}_{k+1}$. 
Note that~\eqref{eq:ol-cov-evolve} describes the open-loop evolution of the measured-state covariance over $[t_k,t_{k+1}]$. Accordingly, the initial condition~\eqref{eq:ol-cov-init} encodes the impulsive update applied at $t_k$ so that this continuous-time segment is consistent with the discrete-time closed-loop covariance recursion in~\eqref{eq:disc-cov-dyn-meas}.

Further, the spacecraft trajectory over $[t_k,t_{k+1}]$ with initial condition $\bm{x}_k$ and control input $\bm{u}^{\mr{cl}}_k$ is given by
\begin{align}
    \check{\bm{x}}_k(t) ={} & \Phi(t,t_k)\bm{x}_k + \Psi(t,t_k)(u_k + K_k(\bm{x}_k+\bm{\zeta}_k-x_k)+\bm{\mu}_k) + \Theta(t,t_k), & & t \in [t_k,t_{k+1}].
\end{align}
Then, for each $t\in[t_k,t_{k+1}]$, the passive-safety chance constraint~\eqref{stoc-ocp:ps-cnstr} imposed on $\check{\bm{x}}_k(t)$ is reformulated as
\begin{align}
& \apskttau^\top \check{x}_k(t) + \bpskttau \ge \cpskttau,~\forall\,\tau\in\tsafespan, \label{ps-chance-reform}
\end{align}
where 
\begin{subequations}
    \begin{align}
        \apskttau \triangleq{} & \rtx{\nabla_1 \sd{\check{\bar{x}}_k(t)}{\BRS}}^\top,\\
        \bpskttau \triangleq{} & \rtx{\sd{\check{\bar{x}}_k(t)}{\BRS}} - \apskttau ^\top \check{\bar{x}}_k(t),\\
        \cpskttau \triangleq{} & \sqrt{Q_{n^x}(\beta^{\mr{ps}})\apskttau^\top\Sigmamk(t)\apskttau}, \label{ps-chance-reform-offset}
    \end{align}\label{ps-chance-reform-parameters}%
\end{subequations}
where satisfaction of~\eqref{ps-chance-reform} implies that
\begin{align}
    \mb{P}(\apskttau^\top \check{\bm{x}}_k(t) + \bpskttau \ge 0,~\forall\,\tau\in\tsafespan) \ge \beta^{\mr{ps}}.
\end{align}
Similarly, for each $t\in[t_k,t_{k+1}]$, the approach-cone chance constraint is reformulated as
\begin{align}
    \aackjt^\top \check{x}_k(t) + \backjt + \cackjt \le 0,~j=1,2, \label{ac-chance-reform}
\end{align}
where
\begin{subequations}
    \begin{align}
        \aackjt \triangleq{} & \rtx{\nabla g^{\mr{ac}}_j(\check{\bar{x}}_k(t))^\top}, \\
        \backjt \triangleq{} & \rtx{g^{\mr{ac}}_j(\check{\bar{x}}_k(t)) - \aackjt^\top \check{\bar{x}}_k(t)}, \\
        \cackjt \triangleq{} & \sqrt{Q_{n^x}(\beta^{\mr{ac}})\aackjt^\top\Sigmamk(t)\aackjt}, \label{ac-chance-reform-offset}
    \end{align}%
\end{subequations}
for $j=1,2$, where~\eqref{ac-chance-reform} implies that
\begin{align}
   \mb{P}(\aackjt^\top \check{\bm{x}}_k(t) + \backjt \le 0,~j=1,2)\ge \beta^{\mr{ac}}.
\end{align}
Note that $Q_{n^x}$ in~\eqref{ps-chance-reform-offset} and~\eqref{ac-chance-reform-offset} is the quantile function of the chi-squared distribution with $n^x$ degrees of freedom.

Given reference states and control inputs, \rtx{the preceding steps} transform~\eqref{stoc-ocp} into the following convex, deterministic optimal control problem in which the nominal control inputs are the decision variables.
\begin{subequations}
\begin{align}
\underset{u_k}{\mr{minimize}}~&~\rtx{\sum_{k=1}^{N-1}}\alpha_k\|u_k\|, & & \label{det-cvx-ocp:obj}\\
\mr{subject~to}~&~x_{k+1} = A_kx_k + B_ku_k + c_k, & & \hphantom{\tau\in\tsafespan,~t\in[t_k,t_{k+1}],~}\,k=1,\ldots,N-1,\label{det-cvx-ocp:dyn}\\
&~\apskttau^\top \check{x}_k(t) + \bpskttau \ge \cpskttau, & & \tau\in\tsafespan,~t\in[t_k,t_{k+1}],~k=1,\ldots,N-1,\label{det-cvx-ocp:ps-cnstr}\\
&~\aackjt^\top \check{x}_k(t) + \backjt + \cackjt \le 0, & & \hspace{0.35cm}j=1,2,~t\in[t_k,t_{k+1}],~k=1,\ldots,N-1,\label{det-cvx-ocp:ac-cnstr}\\
&~(\aactj)^\top u_k \le \bactj - \cactj, & &  \hspace{1.425cm}j=1,\ldots,n^{\mr{act}},~k=1,\ldots,N-1,\label{det-cvx-ocp:ctrl-cnstr}\\
&~\|\selector{r}x_{N_j}\| \le a^+_j,~(e^{\mr{ac}})^\top\selector{r}x_{N_j}\ge a^-_j, & & \hphantom{\tau\in\tsafespan,~t\in[t_k,t_{k+1}],~}\,j=2,3,\label{det-cvx-ocp:phase-transition}\\
&~x_1 = \xinit,~x_N = \xfinal.\label{det-cvx-ocp:bc} 
\end{align}\label{det-cvx-ocp}%
\end{subequations}
%
The state trajectory $\check{x}_k$ in constraints~\eqref{det-cvx-ocp:ps-cnstr} and~\eqref{det-cvx-ocp:ac-cnstr} satisfies~\eqref{eq:lin-ct-dyn} over $[t_k,t_{k+1}]$ with initial condition $x_k$ and control input signal $\nu_k(t,u_k)$, for $t\in[t_k,t_{k+1}]$. \rtx{The chance-constraint reformulation in this section can be further refined by i) explicitly accounting for non-Gaussian state distributions induced by nonlinear dynamics, e.g., with Gaussian mixture models~\cite{boone2022nongaussian}, ii) including a state-estimation filter~\cite{oguri2022stochastic}, and iii) using less conservative deterministic reformulations~\cite{caleb2025chance}.} 

\rtx{Note that LinCov-based methods~\cite{woffinden2022david,goulet2023robust,cavesmith2024angles,cunningham2025robust,margolis2024coopt,margolis2024robust} also rely on linearized propagation of the state mean and covariance, as in~\eqref{eq:lin-disc-dyn} and~\eqref{eq:disc-cov-dyn-meas}. In the gradient-free formulations of~\cite{woffinden2022david,goulet2023robust,cavesmith2024angles,cunningham2025robust}, this covariance propagation is performed for each candidate trajectory evaluated by the underlying particle-swarm or genetic-algorithm search, which can lead to many covariance propagations per iteration. The gradient-based formulations in~\cite{margolis2024coopt,margolis2024robust} are closer in spirit to our approach in that they also embed covariance propagation within the optimization and exploit gradient information. Our method differs in how the propagated covariance is used. In the cited formulations, LinCov is used primarily to define covariance-based performance and safety measures, often expressed as 3$\sigma$ bounds or margins, whereas we use the propagated covariance to construct ellipsoidal confidence set-based tightening \cite[Rem. 1]{lew2020chance} of the original constraints at each SCP iteration, thereby approximating the chance constraints.} Although all constraints in~\eqref{det-cvx-ocp} are convex, two aspects make the problem difficult to solve directly: i) the signed distance to the BRS of a time-varying nonlinear system is required in~\eqref{det-cvx-ocp:ps-cnstr}, and ii) the path constraints~\eqref{det-cvx-ocp:ps-cnstr} and~\eqref{det-cvx-ocp:ac-cnstr} must hold at infinitely many time instants. We address these challenges in the following subsections.
%
\subsection{Approximation of BRS}\label{subsec:approx-BRS}
\rtx{The reformulated passive-safety constraint~\eqref{det-cvx-ocp:ps-cnstr}, together with the definitions in~\eqref{ps-chance-reform-parameters}, requires the signed distance to the BRS of the time-varying nonlinear system~\eqref{visit-sc-dyn}. In general, computing BRS is intractable for nonlinear systems with state dimension greater than four~\cite[Fig. 2]{chen2018hamiltonjacobi}. For affine dynamical systems, however, BRS can be computed efficiently by numerical integration of first-order sensitivities. We therefore approximate the BRS for~\eqref{visit-sc-dyn} by linearizing the free-drift motion about a reference state trajectory. We also use polyhedral over-approximations of the avoid sets in Figure~\ref{fig:nrho-rdv-formulation}, because the signed-distance function and its gradient can be computed efficiently for such sets (see Appendix~\ref{app:signed-distance}). Thus, $\mc{A}$ is represented by a polyhedral open set $\{z\in\bR^{n^x}\,|\,\Havoid z < \havoid\}$.}

Let $x_1,\dots,x_N$ and $u_1,\ldots,u_{N-1}$ denote nominal states and control inputs, respectively. For each $k=1,\ldots,N-1$, let state trajectory $\check{x}_{k}$ satisfy
\begin{subequations}
\begin{align}
\dot{\check{x}}{}_k(t) ={} &f(t,\check{x}_k(t),\nu_k(t,u_k)), & & t\in[t_k,t_{k+1}],\\
\check{x}_k(t_k) ={} & x_k.
\end{align}\label{dt2ct-traj}%
\end{subequations}
The free-drift trajectory \rtx{starting at $t\in [t_k,t_{k+1}]$ from $\check{x}_k(t)$, denoted by $\check{x}_{k,t}$,} satisfies
\begin{subequations}
\begin{align}
    \derv{\check{x}}{}_{k,t}(\tau) ={} & f(t+\tau,\check{x}_{k,t}(\tau),\zeros{n^u}), & & \tau\in\tsafespan,\label{free-drift-dyn-t:ode}\\
    \check{x}_{k,t}(0) ={} & \check{x}_k(t).
\end{align}\label{free-drift-dyn-t}%
\end{subequations}
Then, the passive-safety constraint~\eqref{ps-cnstr-BRS-naive} can be equivalently expressed as
\begin{align}
    \check{x}_{k,t}(\tau) \notin \mc{A},\quad\forall\,\tau\in\tsafespan,~t\in[t_k,t_{k+1}],~k=1,\ldots,N-1.\label{ps-cnstr-free-drift-k}
\end{align}
\rtx{Given reference states $\bar{x}_1,\ldots,\bar{x}_N$ and control inputs $\bar{u}_1,\ldots,\bar{u}_{N-1}$, we define the reference state trajectories $\check{\bar{x}}_k$ and $\check{\bar{x}}_{k,t}$ using~\eqref{dt2ct-traj} and~\eqref{free-drift-dyn-t}, respectively. The free-drift motion is then approximated by linearizing~\eqref{free-drift-dyn-t:ode} about the reference free-drift trajectory $\check{\bar{x}}_{k,t}$}
\begin{subequations}
\begin{align}
    {\derv{\check{x}}}_{k,t}(\tau) ={} & f(t+\tau,\check{x}_{k,t}(\tau),\zeros{n^u}),\\
    \approx{} & f(t+\tau,\check{\bar{x}}_{k,t}(\tau),\zeros{n^u}) + \nabla_2 f(t+\tau,\check{\bar{x}}_{k,t}(\tau),\zeros{n^u}) (\check{{x}}_{k,t}(\tau)-\check{\bar{x}}_{k,t}(\tau)),\\
    ={} & \check{A}_{k,t}(\tau)\check{{x}}_{k,t}(\tau) + \check{c}_{k,t}(\tau),\label{free-drift-linsys}
\end{align}\label{approx-free-drift}%
\end{subequations}
where $\check{A}_{k,t}(\tau) \triangleq \nabla_2 f(t+\tau,\check{\bar{x}}_{k,t}(\tau),\zeros{n^u})$ and $\check{c}_{k,t}(\tau)\triangleq f(t+\tau,\check{\bar{x}}_{k,t}(\tau),\zeros{n^u}) - \check{A}_{k,t}(\tau)\check{\bar{x}}_{k,t}(\tau)$. Next, the first-order sensitivities of~\eqref{free-drift-dyn-t:ode} over the safety horizon can be computed by solving
\begin{subequations}
\begin{align}
{\derv{\check{\bar{x}}}}_{k,t}(\tau) ={} & f(t+\tau,\check{\bar{x}}_{k,t}(\tau),\zeros{n^u}), & & \tau\in\tsafespan,\\
{\derv{\Phi}}_{k,t}(\tau) ={} & \check{A}_{k,t}(\tau){\Phi}_{k,t}(\tau), & & \tau \in \tsafespan,\\
{\derv{\Theta}}_{k,t}(\tau) ={} & \check{A}_{k,t}(\tau){\Theta}_{k,t}(\tau) + \check{c}_{k,t}(\tau), & & \tau \in \tsafespan,\\  
{\check{\bar{x}}}_{k,t}(0) ={} & \check{\bar{x}}_k(t),\\
{\Phi}_{k,t}(0) ={} & \eye{n^x},\\
{\Theta}_{k,t}(0) ={} & \zeros{n^x}.
\end{align}
\end{subequations}
The free-drift state $\check{x}_{k,t}(\tau)$ is approximated as
\begin{align}
    \check{x}_{k,t}(\tau) \approx {\Phi}_{k,t}(\tau)\check{x}_k(t) + {\Theta}_{k,t}(\tau).\label{eq:free-drift-state-approx}
\end{align}
Then, the passive-safety constraint~\eqref{ps-cnstr-free-drift-k} is approximated as
\begin{align}   
    {\Phi}_{k,t}(\tau)\check{x}_{k}(t) + {\Theta}_{k,t}(\tau) \notin \mc{A},\quad\forall\,\tau\in\tsafespan,~t\in[t_k,t_{k+1}],~k=1,\ldots,N-1,
\end{align}
which is equivalent to
\begin{align}
    \check{x}_k(t) \notin \approxBRS,\quad \forall \tau\in\tsafespan,~t\in [t_k,t_{k+1}],~k=1,\ldots,N-1,\label{ps-cnstr-approx-BRS-naive}
\end{align}
where
\begin{align}
\approxBRS \triangleq \{z\,|\,\Havoid\,{\Phi}_{k,t}(\tau)z<\havoid-\Havoid\,{\Theta}_{k,t}(\tau)\},
\end{align}
is an approximation of the BRS $\BRS$, which we use in place of $\BRS$ in $\apskttau$, $\bpskttau$, and $\cpskttau$. Note that $\approxBRS$ is an open set, since $\mc{A}$ is open.
%
\subsection{Isoperimetric Reformulation of Path Constraints}\label{subsec:isoperi}
To ensure that the path constraints are satisfied in continuous time, we use an isoperimetric reformulation~\cite[Sec. 10]{hartl1995survey}. Specifically, for each $k = 1,\ldots,N-1$, the pointwise-in-time path constraints~\eqref{det-cvx-ocp:ps-cnstr} and~\eqref{det-cvx-ocp:ac-cnstr} can be equivalently enforced using integral equality constraints
\begin{subequations}
\begin{align}
    & \int_{t_k}^{t_{k+1}}\int_{0}^{\tsafe}\big|\!-\!\apskttau^\top \check{x}_k(t) - \bpskttau + \cpskttau\big|_+^2\mr{d}\tau\,\mr{d}t = 0,\label{eq:isoperi-path-cnstr:ps}\\
    & \int_{t_k}^{t_{k+1}}\big|\aackjt^\top \check{x}_k(t)+\backjt+\cackjt\big|_+^2\mr{d}t = 0, & & j =1,2, 
\end{align}\label{eq:isoperi-path-cnstr}%
\end{subequations}
where each integrand is the composition of a path constraint function and the differentiable exterior penalty function $\square \mapsto |\square|_+^2$. \rtx{Each integrand is continuous, nonnegative, and equal to zero if and only if the corresponding path constraint is satisfied. Therefore, the integral equals zero if and only if the integrand equals zero for all times in the interval (except on sets of measure zero). Thus,} satisfying~\eqref{eq:isoperi-path-cnstr} is equivalent to satisfying~\eqref{det-cvx-ocp:ps-cnstr} and~\eqref{det-cvx-ocp:ac-cnstr} for all $t\in[t_k,t_{k+1}]$.

In practice, the integral constraints in~\eqref{eq:isoperi-path-cnstr} can be enforced \rtx{as a} boundary value problem with an auxiliary system
\begin{subequations}
    \begin{align}
        \dot{\check{y}}{}_k(t) ={} & \begin{bmatrix} \displaystyle\int_{0}^{\tsafe}\big|\!-\!\apskttau^\top \check{x}_k(t) - \bpskttau + \cpskttau\big|_+^2\mr{d}\tau \\ \big|\aackt{1}^\top \check{x}_k(t)+\backt{1}+\cackt{1}\big|_+^2\\[0.2cm]
            \big|\aackt{2}^\top \check{x}_k(t)+\backt{2}+\cackt{2}\big|_+^2 \end{bmatrix},\label{eq:cnstr-integrator} & & t\in[t_k,t_{k+1}],\\
        \check{y}_k(t_{k+1}) ={} & \check{y}_k(t_k). \label{eq:cnstr-integrator-bc} 
    \end{align}\label{eq:aux-cnstr-integrator-bvp}%
\end{subequations}
\rtx{For further discussion of the equivalence among~\eqref{det-cvx-ocp:ps-cnstr}--\eqref{det-cvx-ocp:ac-cnstr},~\eqref{eq:isoperi-path-cnstr}, and~\eqref{eq:aux-cnstr-integrator-bvp}}, we refer to~\cite[Lem. 2]{elango2025ctscvx}.

The linearization of~\eqref{eq:cnstr-integrator} with respect to the reference states and control inputs, along with discretization over $t_1,\ldots,t_N$, yields
\begin{align}
    y_{k+1} = y_k + G_k^x x_k + G_k^uu_k + h_k,
\end{align}
%
%
where
\begin{subequations}
\begin{align}
    G^x_k \triangleq{} & \int_{t_k}^{t_{k+1}} \check{G}_k(t)\Phi(t,t_k)\mr{d}t,\\
    G^u_k \triangleq{} & \int_{t_k}^{t_{k+1}} \check{G}_k(t)\Psi(t,t_k)\mr{d}t,\\
    h_k \triangleq{} & \int_{t_k}^{t_{k+1}}  \check{h}_k(t)\mr{d}t,\\  
    \check{G}_k(t) \triangleq{} &    \begin{bmatrix}
                        -2\displaystyle\int_0^{\tsafe}\big|\!-\!\rtx{\sd{\check{\bar{x}}_k(t)}{\approxBRS}} + \cpskttau\big|_+\apskttau^\top\mr{d}\tau\\
                        2\big|g^{\mr{ac}}_1(\check{\bar{x}}_k(t))+\cackt{1}\big|_+\aackt{1}^\top\\[0.2cm]
                        2\big|g^{\mr{ac}}_2(\check{\bar{x}}_k(t))+\cackt{2}\big|_+\aackt{2}^\top
                    \end{bmatrix},\\
    \check{h}_k(t) \triangleq{} & \begin{bmatrix}
        \displaystyle\int_{0}^{\tsafe}\big|\!-\!\rtx{\sd{\check{\bar{x}}_k(t)}{\approxBRS}} + \cpskttau\big|_+^2 + 2\big|\!-\!\rtx{\sd{\check{\bar{x}}_k(t)}{\approxBRS}} + \cpskttau\big|_+\apskttau^\top\check{\bar{x}}_k(t)\mr{d}\tau\\
        \big| g^{\mr{ac}}_1(\check{\bar{x}}_k(t)) + \cackt{1}\big|_+^2 - 2\big|\rtx{g^{\mr{ac}}_1}(\check{\bar{x}}_k(t))+\cackt{1}\big|_+\aackt{1}^\top\check{\bar{x}}_k(t)\\[0.2cm]
        \big| g^{\mr{ac}}_2(\check{\bar{x}}_k(t)) + \cackt{2}\big|_+^2 - 2\big|\rtx{g^{\mr{ac}}_2}(\check{\bar{x}}_k(t))+\cackt{2}\big|_+\aackt{2}^\top\check{\bar{x}}_k(t)
                 \end{bmatrix}.\label{ct-linearize-G}
\end{align}
\end{subequations}
When using ZOH or FBP control inputs, enforcing~\eqref{eq:isoperi-path-cnstr:ps} guarantees passive safety in the event of an underburn (i.e., the thruster fires for a shorter duration than intended). Consider FBP. Under an underburn, the actual burn ends at $t_k+t'$ with $t'\in[0,\tburn)$. Since the integrand in~\eqref{eq:isoperi-path-cnstr:ps} is nonnegative and
\begin{align}
\int_{t_k}^{t_k+\tburn}\!\!\int_{0}^{\tsafe}
\big|\!-\!\apskttau^\top \check{x}_k(t) - \bpskttau + \cpskttau\big|_+^2
\,\mr{d}\tau\,\mr{d}t \;=\; 0,
\end{align}
it follows that the integrand is zero almost everywhere on $[t_k,\,t_k+\tburn]\times[0,\tsafe]$, and hence, for every $t'\in[0,\tburn]$,
\begin{align}
\int_{t_k}^{t_k+t'}\!\!\int_{0}^{\tsafe}
\big|\!-\!\apskttau^\top \check{x}_k(t) - \bpskttau + \cpskttau\big|_+^2
\,\mr{d}\tau\,\mr{d}t \;=\; 0.
\end{align}
Therefore, if the thruster shuts off at any time in $[t_k,t_k+\tburn]$, the resulting free-drift trajectory remains passively safe. In contrast, under a purely impulsive actuation model, an ``underburn'' (i.e., \rtx{delivery of} only a fraction of the commanded impulse) can yield unsafe free drift. Here, underburn for FBP and ZOH \rtx{refers specifically} to the thruster firing for less than the commanded duration. Cases involving incorrect thrust orientation are not considered in this work.
%
\subsection{Sequential Convex Programming}
\begin{figure}[!ht]
\centering
\includegraphics[width=0.85\linewidth]{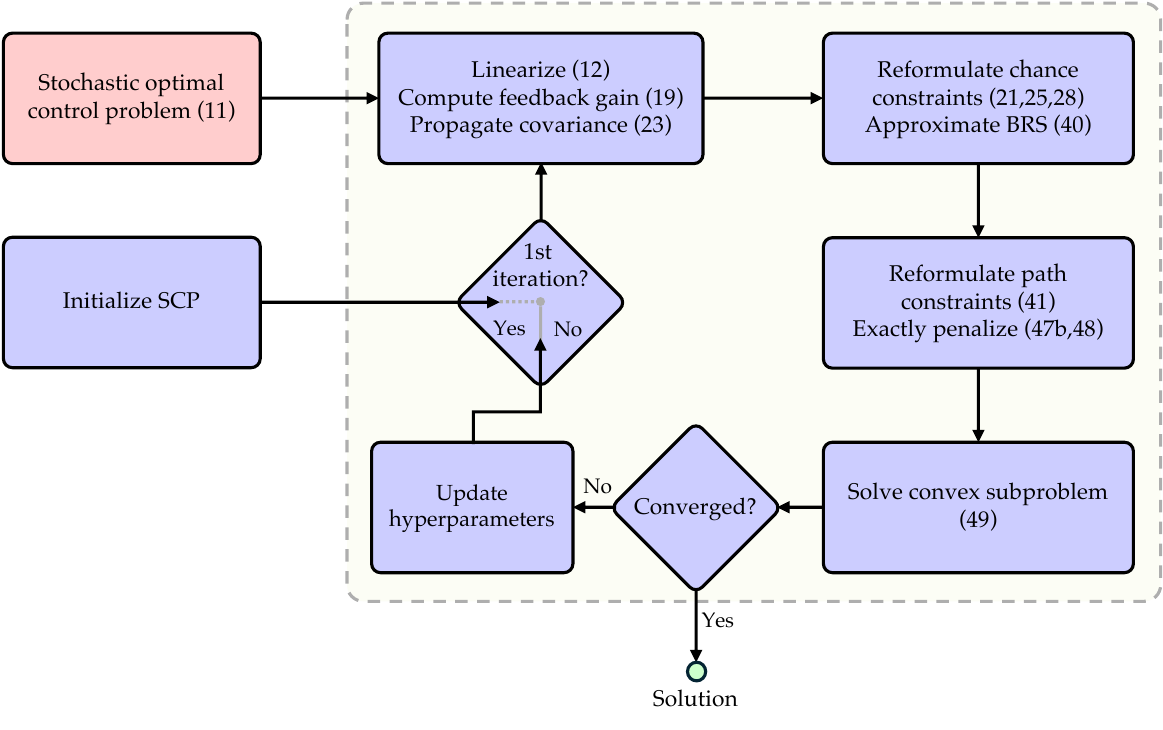}    
\caption{SCP-based approach for computing the rendezvous maneuver in problem~\eqref{stoc-ocp}.}\label{fig:scp-blkdiag}%
\end{figure}

\rtx{We solve~\eqref{stoc-ocp} using an SCP-based approach. First, after approximating the BRS and applying the isoperimetric reformulation of the path constraints, we obtain the following finite-dimensional convexified optimal control problem}
\begin{subequations}
    \begin{align}
    \underset{u_k}{\mr{minimize}}~&~\sum_{k=1}^{N-1} \alpha_k\|u_k\|, & & \\
    \mr{subject~to}~&~x_{k+1} = A_kx_k + B_ku_k + c_k, & & \hphantom{j = 1,\ldots,n^{\mr{act}},~}\,k = 1,\ldots,N-1,\label{det-cvx-isoperi-ocp:dyn}\\
    ~&~G^x_kx_k + G^u_ku_k + h_k = \zeros{3}, & & \hphantom{j = 1,\ldots,n^{\mr{act}},~}\,k = 1,\ldots,N-1,\label{det-cvx-isoperi-ocp:ac-ps-cnstr}\\
    ~&~(\aactj)^\top u_k \le \bactj - \cactj, & & j = 1,\ldots,n^{\mr{act}},~k = 1,\ldots,N-1,\label{det-cvx-isoperi-ocp:ctrl-cnstr}\\
    &~\|\selector{r}x_{N_j}\| \le a^+_j,~(e^{\mr{ac}})^\top\selector{r}x_{N_j}\ge a^-_j, & & \hphantom{\tau\in\tsafespan,~t\in[t_k,t_{k+1}],~}\,j=2,3,\label{det-cvx-isoperi-ocp:phase-transition}\\
    &~x_1 = \xinit,~x_N = \xfinal.\label{det-cvx-isoperi-ocp:bc}        
    \end{align}\label{det-cvx-isoperi-ocp}%
\end{subequations}
Next, we relax~\eqref{det-cvx-isoperi-ocp:ac-ps-cnstr} to a linear inequality with a positive tolerance $\epsilon$
\begin{align}
    G^x_kx_k + G^u_ku_k + h_k \le \epsilon\ones{3}, & & k = 1, \ldots, N-1,\label{ac-ps-cnstr-relax}
\end{align}
to avoid automatic violation of constraint qualifications at feasible points of~\eqref{det-cvx-isoperi-ocp} (see~\cite[Sec. 3]{elango2025ctscvx} for further details). However, the relaxation in~\eqref{ac-ps-cnstr-relax} can admit solutions that are infeasible with respect to the path constraints. We mitigate this \rtx{effect} by systematically tightening the path constraints using the approach in~\cite[Thm.~10]{elango2025ctscvx}, which counteracts the relaxation with only mild conservatism. The amount of tightening increases with the cube root of the small tolerance $\epsilon$. \rtx{Finally, we replace~\eqref{det-cvx-isoperi-ocp:dyn} and~\eqref{ac-ps-cnstr-relax} with exact penalties using the $\ell_1$ penalty functions $\square\mapsto\|\square\|_1$ and $\square\mapsto|\square|_+$. These penalties avoid artificial infeasibility caused by linearizing the nonlinear spacecraft dynamics and the nonconvex path constraints~\cite[Sec. III.A.4]{szmuk2020successive}. We also include a quadratic proximal term in the objective function to limit the step size of the current SCP iterate. These steps result in the following convex subproblem}
\begin{subequations}
\begin{align}
\underset{u_k}{\mr{minimize}}~&~\sum_{k=1}^{N-1}(\alpha_k\|u_k\| + w^{\mr{ep}}\ones{n^x}^\top(\omega^+_k + \omega^-_k) + w^{\mr{ep}}\ones{3}^\top \rho_k) & & \\
&~+w^{\mr{px}}\|x_N-\bar{x}_N\|^2+w^{\mr{px}}\sum_{k=1}^{N-1}(\|x_k-\bar{x}_k\|^2 + \|u_k-\bar{u}_k\|^2)\nonumber\\
\mr{subject~to}~&~x_{k+1} = A_kx_k + B_ku_k + c_k + \omega^+_k - \omega^-_k, & & \hphantom{j = 1,\ldots,n^{\mr{act}},~}\,k = 1,\ldots,N-1,\label{det-cvx-isoperi-rlx-ocp:dyn}\\
~&~G^x_kx_k + G^u_ku_k + h_k \le \epsilon\ones{3} + \rho_k, & & \hphantom{j = 1,\ldots,n^{\mr{act}},~}\,k = 1,\ldots,N-1,\label{det-cvx-isoperi-rlx-ocp:ac-ps-cnstr}\\
~&~\omega_k^+ \ge \zeros{n^x},~\omega^-_k \ge \zeros{n^x},~\rho_k\ge\zeros{3}, & & \hphantom{j = 1,\ldots,n^{\mr{act}},~}\,k = 1,\ldots,N-1,\\
~&~(\aactj)^\top u_k \le \bactj - \cactj, & & j = 1,\ldots,n^{\mr{act}},~k = 1,\ldots,N-1,\label{det-cvx-isoperi-rlx-ocp:ctrl-cnstr}\\
&~\|\selector{r}x_{N_j}\| \le a^+_j,~(e^{\mr{ac}})^\top\selector{r}x_{N_j}\ge a^-_j, & & \hphantom{\tau\in\tsafespan,~t\in[t_k,t_{k+1}],~}\,j=2,3,\label{det-cvx-isoperi-rlx-ocp:phase-transition}\\
&~x_1 = \xinit,~x_N = \xfinal,\label{det-cvx-isoperi-rlx-ocp:bc}        
\end{align}\label{det-cvx-isoperi-rlx-ocp}%
\end{subequations}
which can be \rtx{solved efficiently} with an off-the-shelf conic optimization solver~\cite{domahidi2013ecos,mosek,yu2022proportionalintegral}. The weights for the $\ell_1$-penalty slack variables and the proximal term are $w^{\mr{ep}}$ and $w^{\mr{px}}$, respectively. \rtx{The previous SCP iterate, $\bar{x}_1,\dots,\bar{x}_N$ and $\bar{u}_1,\dots,\bar{u}_{N-1}$, is used for the linearization steps in Sections~\ref{subsec:cov-evolve-feedback}--\ref{subsec:isoperi}. Given positive tolerances $\epsilon^{\mr{ep}}$ and $\epsilon^{\mr{px}}$, SCP terminates when}
\begin{subequations}
\begin{align}
    & \sum_{k=1}^{N-1} \rtx{\big(\ones{n^x}^\top(\omega^+_k+\omega^-_k) + \ones{3}^\top\rho_k\big)} \le \epsilon^{\mr{ep}},\\
    & \|x_N-\bar{x}_N\| + \sum_{k=1}^{N-1}\rtx{\big(\|x_k-\bar{x}_k\|^2 + \|u_k-\bar{u}_k\|^2\big)}  \le \epsilon^{\mr{px}}.
\end{align}\label{term-crit}%
\end{subequations}
For appropriate choices of $w^{\mr{px}}$, using a proximal term is equivalent to imposing a hard trust-region constraint with a specific trust-region radius~\cite{cartis2011evaluation}. The complete SCP approach is summarized in Figure~\ref{fig:scp-blkdiag}. The iterative part of the algorithm is enclosed in the dashed box, and the ``Update hyperparameters'' block denotes the use of methods for tuning $w^{\mr{px}}$ and $w^{\mr{ep}}$ to improve practical SCP convergence~\cite{cartis2011evaluation,xie2025high,mceowen2025autotuned}.

The proposed SCP-based approach is not guaranteed to find a feasible solution from every initial condition. However, convergence can be verified before the rendezvous maneuver starts. When the method converges, the resulting solution satisfies the path constraints (passive safety and the approach cone) and \rtx{remains safe under underburns for} the ZOH or FBP parametrizations.

\subsection{Passive-Safety at SCP Convergence}
We now quantify passive-safety constraint satisfaction at SCP convergence. Let
\begin{align*}
    Z ={} & (x_1,\ldots,x_N,\,u_1,\ldots,u_{N-1})\in\bR^{n^Z},\\
    Y ={} & (\omega^+_1,\omega^-_1,\ldots,\omega^+_{N-1},\omega^-_{N-1},\rho_1,\ldots,\rho_{N-1})\in\bR^{n^Y},
\end{align*}
collect the states and control inputs, and the slack variables, respectively, of the convex subproblem~\eqref{det-cvx-isoperi-rlx-ocp}. Let $\bar{Z}$ denote the stacked state/control input solution from the previous SCP iterate. At SCP convergence (i.e., as $Z\to\bar{Z}$ and $Y\to\zeros{n^Y}$), \rtx{the free-drift approximation} in~\eqref{eq:free-drift-state-approx} becomes exact. Hence, for each $k=1,\ldots,N-1$, for all $t\in[t_k,t_{k+1}]$ and $\tau\in\tsafespan$,
\begin{align}
    \Phi_{k,t}(\tau)\,\check{x}_k(t) + \Theta_{k,t}(\tau) \notin \mc{A}
    \;\;\implies\;\;
    \check{x}_{k,t}(\tau)\notin\mc{A},    
\end{align}
which is equivalent to
\begin{align}
    \check{x}_k(t)\notin \approxBRS \;\;\implies\;\; \check{x}_k(t)\notin \BRS.
    \label{eq:sd-imply-at-SCP-conv}
\end{align}
Note that~\eqref{eq:sd-imply-at-SCP-conv} does not assert $\approxBRS=\BRS$ at convergence. It is a property of the converged trajectory, not a set equality.
%

Moreover, as $Z\to\bar{Z}$ the affine approximation of the signed-distance of the free-drift state from the avoid set becomes exact
\begin{align}
    \apskttau^\top \check{\bar{x}}_k(t) + \bpskttau \;=\; \rtx{\sd{\check{\bar{x}}_k(t)}{\approxBRS}},
    \quad \forall\,t\in[t_k,t_{k+1}],\,k=1,\ldots,N-1.
    \label{eq:sd-term-at-SCP-conv}
\end{align}
Combining~\eqref{eq:sd-imply-at-SCP-conv} and~\eqref{eq:sd-term-at-SCP-conv}, the following holds at convergence:
\begin{align}
    \rtx{\sd{\check{\bar{x}}_k(t)}{\BRS}} \;\ge\; \cpskttau,
    \quad \forall\,\tau\in\tsafespan,\,t\in[t_k,t_{k+1}],\,k=1,\ldots,N-1.
    \label{eq:ps-cnstr-SCP-conv}
\end{align}
Thus, the converged nominal trajectory is passively safe with respect to the true BRS of the avoid set, with margin $\cpskttau$ accounting for uncertainty. 

In highly nonlinear regions of the NRHO, such as near the perilune, errors from linearizing the free-drift motion can be significant, potentially leading to slower SCP convergence. This limitation does not affect the apolune rendezvous scenario considered here.

%
\section{Numerical Results}\label{sec:num-results}
This section demonstrates the proposed SCP-based approach through a numerical case study of rendezvous to the Gateway, based on the specifications in Section~\ref{sec:prb-form}. \rtx{Table~\ref{tab:phases-data} lists the case study parameters, including} i) the covariances of the initial state, relative-range measurements at the decision points, and actuation error, ii) chance-constraint confidence levels, iii) \rtx{approach-cone parameters and range constraints} at D\textsubscript{2} and D\textsubscript{3}, and iv) the initial and final conditions of the nominal spacecraft state. These values are based on related case studies in~\cite{nakamura2023rendevous}. The rotation matrix $\rotmatLVLH$ transforms position and velocity from the Sun Referenced LVLH frame~\cite{irsis2019} to the Gateway-centered J2000 inertial frame. The polytopic avoid set in phase $j$, given by
\begin{table}[!htp]
\centering
\caption{Parameter values for the three-phase rendezvous maneuver to the Gateway}\label{tab:phases-data}
\setstretch{1.5}
\begin{tabular}{lc}
    \hline\hline
    \textbf{Parameter} & \textbf{Value}\\\hline
    $N_2$, $N_3$, $N$ & $4$, $8$, $12$\\
    $t_1,\ldots,t_N$ [hr] & $0$,~$30$,~$38$,~$42$,~$43$,~$44$,~$45$,~$46$,~$47.25$,~$47.5$,~$47.75$,~$48$\\
    $\tsafe$ [hr] & $24$\\
    $\tburn$ [min] & $15$ (phases 1 \& 2), $7.5$ (phase 3)\\
    $r^{\scriptscriptstyle\mc{A}_1}$, $r^{\scriptscriptstyle\mc{A}_2}$, $r^{\scriptscriptstyle\mc{A}_3}$ [km]& $10$, $1$, $0.2$\\
    $\Sigma^{\mr{i}}$ [km$^2$,\,km$^2$/hr$^2$] & $\mr{diag}(33.33,33.33,33.33,6,6,6)^2$\\
    $\Sigma^{\mr{rr}}_{1}$ [km$^2$,\,km$^2$/hr$^2$] & $\mr{diag}(6.66,6.66,6.66,0.25,0.25,0.25)^2$\\
    $\Sigma^{\mr{rr}}_{N_2}$ [km$^2$,\,km$^2$/hr$^2$] & $\mr{diag}(0.08,0.08,0.08,0.04,0.04,0.04)^2$\\
    $\Sigma^{\mr{rr}}_{N_3}$ [km$^2$,\,km$^2$/hr$^2$] & $\mr{diag}(0.009,0.009,0.009,0.03,0.03,0.03)^2$\\
    $\Sigma^{\mr{rr}}_{N}$ [km$^2$,\,km$^2$/hr$^2$] & $\mr{diag}(0.006,0.006,0.006,0.01,0.01,0.01)^2$\\
    $\Sigma^{\mr{act}}$ [km$^2$/hr$^4$] & $\mr{diag}(0.5,0.5,0.5)^2$\\
    $\beta^{\mr{ps}},~\beta^{\mr{ac}}$ & $0.8$, $0.8$\\
    $\theta^{\mr{ac}},~e^{\mr{ac}}$ & $55^\circ$, $\rotmatLVLH\,\selector{r}^\top(0,0,1)$\\
    $u^{\max}$ [km/hr$^2$] & 120 \\
    $\xinit$, $\xfinal$ [km,\,km/hr]& $\rotmatLVLH(0, -600, 800,2.5,30,-20)$, $\rotmatLVLH(0,0,0.5,0,0,0)$\\
    $a^+_2$, $a^-_2$, $a^+_3$, $a^-_3$ [km] & $55$, $45$, $6.5$, $3.5$\\\hline\hline
\end{tabular}
\end{table}
\begin{align}
    \mc{A}_j = \left\{ x\in\bR^{n^x}\,\middle|\,\left[\begin{bmatrix} \eye{3} \\ -\eye{3} \end{bmatrix}\selector{r}\rotmatLVLH\selector{r}^\top~~~\zeros{n^x\times 3}\right]x < r^{\scriptscriptstyle\mc{A}_j}\ones{n^x} \right\},
\end{align}
for $j=1,2,3$, forms an open box in the position coordinates, where $r^{\scriptscriptstyle\mc{A}_j}$ is the half edge length of the box. \rtx{We use FBP control inputs to model realistic pulsed thrusters, and choose the admissible control set as an infinity-norm ball, i.e., the polytope~$\mc{U} = \{u\in\bR^3\,|\,\|u\|_\infty \le u^{\max}\}$.}

Navigation accuracy depends on the geometry of the rendezvous trajectory (e.g., relative range and sensor/boresight orientation) and the sensing modality (e.g., ground-based tracking vs.\ Gateway-relative vision/LiDAR). This coupling introduces a challenging nonconvexity in the guidance problem because the covariance depends on the evolving spacecraft state. At far range (near D\textsubscript{1}), the spacecraft relies on ground-based orbit determination; at mid range (beyond D\textsubscript{2}), it switches to line-of-sight (LOS) and Gateway-relative measurements. Accordingly, we specify the covariance of the relative-range measurements at the decision points and linearly interpolate the boundary values between successive decision points. This simplification is physically reasonable because the spacecraft's range decreases roughly monotonically from D\textsubscript{1} to D\textsubscript{4}, and it preserves tractability of the SCP and the continuous-time chance-constraint representation. Furthermore, in this case study, the spacecraft approaches along the Sun-Gateway axis, maintaining target illumination and partially mitigating vision-sensor variability. The proposed framework is compatible with any time-varying covariance profile. More detailed, state-dependent, and multimodal navigation-error models can be incorporated through an outer loop (e.g., freezing the covariances from the previous SCP iterate) or via a scenario-based approach. A natural extension is to embed a navigation filter (e.g., an extended Kalman filter) within the framework to capture state-dependent navigation errors more accurately while avoiding excessive conservatism.

\rtx{Operational requirements also constrain the frequency of control inputs, particularly in phase~1~\cite{nakamura2023rendevous}. Reducing the frequency of control actions is desirable in the initial phases, when measured state dispersions are large, because frequent correction maneuvers can lead to large-magnitude closed-loop control inputs. Reducing the total number of control actions also lowers the burden on ground operators and the likelihood of thruster failure. We therefore selected a time-discretization grid with the minimum number of nodes that satisfied these guidelines. The correction maneuvers coincide with scheduled nominal control inputs rather than being implemented as separate additional burns.}

Because of challenging nonconvex constraints, such as the continuous-time passive-safety constraint, the initial guess significantly affects SCP convergence. We initialize the state with a straight-line profile (linear interpolation between the nominal states at D\textsubscript{1} and D\textsubscript{4} across intermediate nodes) and set the control input to zero on an initially denser time grid (20 nodes in the case study). After obtaining a feasible solution, we gradually prune the grid \rtx{according to the operational guidelines above, then warm-start SCP on the sparser grid and repeat this process until convergence. This warm-start procedure improves convergence to physically sensible solutions.} 

Two layers of scaling are used within the SCP algorithm to ensure reliable numerical performance. i) Scale--then--discretize: the continuous-time dynamics model~\eqref{cislunar-sc-dyn} is scaled using physically meaningful units (hr, km, km/hr) so that state and control input profiles have similar orders of magnitude, improving integration and linearization stability~\cite{ross2018scaling}. ii) Discretize--then--scale: after constructing the convex subproblem~\eqref{det-cvx-isoperi-rlx-ocp}, we apply affine scaling to the decision variables to obtain comparable magnitudes and precondition the constraint matrices using row normalization~\cite{kamath2025optimal}, which improves SCP convergence~\cite{sagliano2014performance}. After scaling and preconditioning, choosing $\epsilon$ in~\eqref{ac-ps-cnstr-relax} on the order of $10^{-3}$--$10^{-4}$ yields no physically meaningful violation of the original constraint~\eqref{det-cvx-isoperi-ocp:ac-ps-cnstr}.

We performed a Monte Carlo simulation with 1000 samples using a nominal three-phase maneuver computed by the proposed method. \rtx{Table~\ref{tab:mc} reports the fuel-consumption range} and the degree of constraint satisfaction observed across the samples. Notably, \rtx{no sample violates} the passive-safety constraint; all \rtx{remain safe even under underburns}. This outcome reflects the conservativeness of the chance-constrained formulation described in Section~\ref{subsec:chance-cnstr-reform}: even relatively small values of $\beta^{\mr{ps}}$ (e.g., around 0.8) can yield a high percentage of safe outcomes in Monte Carlo simulations. For a detailed discussion of the conservativeness of chance-constraint formulations, we refer the reader to~\cite{lew2020chance}. For some initial conditions, no feasible rendezvous maneuver may exist. The initial condition in Table~\ref{tab:phases-data} is selected to ensure feasibility (including underburn safety) for the nominal case and for an anticipated dispersion of initial conditions about it, as would be done during mission design for the rendezvous.

\begin{figure}[!htp]
\centering
\begin{subfigure}[b]{0.49\linewidth}
\centering
\includegraphics[width=0.89\linewidth]{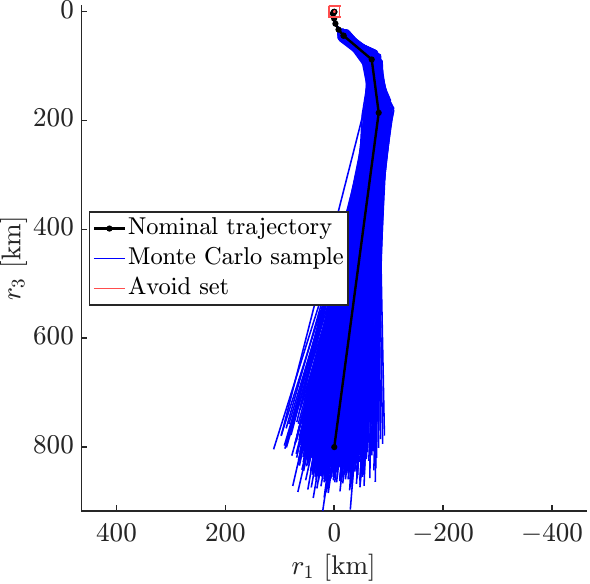}
\end{subfigure}
\begin{subfigure}[b]{0.49\linewidth}
\centering
\includegraphics[width=0.89\linewidth]{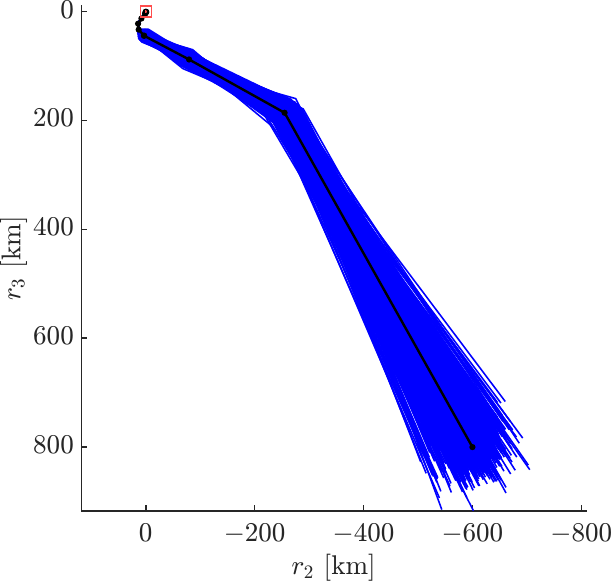}
\end{subfigure}

\begin{subfigure}[b]{0.49\linewidth}
\centering
\includegraphics[width=0.85\linewidth]{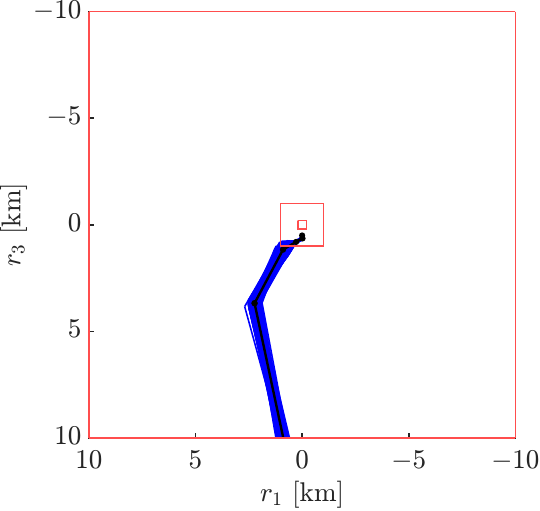}
\end{subfigure}
\begin{subfigure}[b]{0.49\linewidth}
\centering
\includegraphics[width=0.85\linewidth]{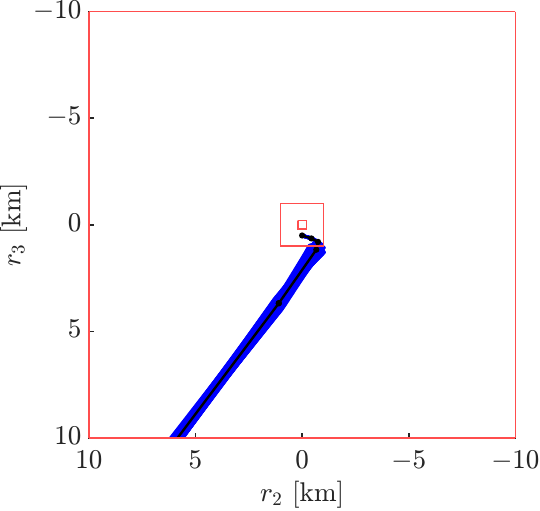}
\end{subfigure}

\begin{subfigure}[b]{0.49\linewidth}
\centering
\includegraphics[width=0.85\linewidth]{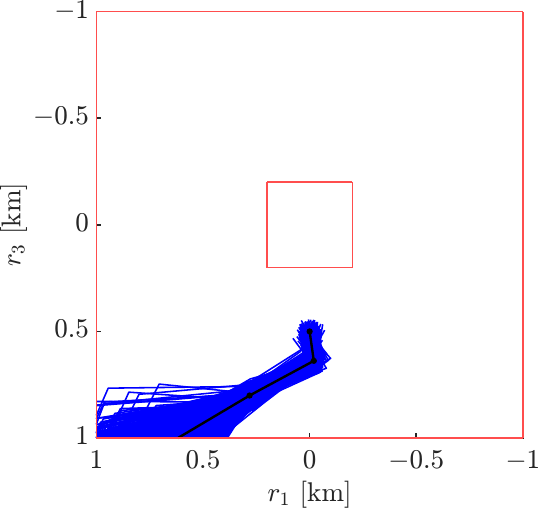}
\end{subfigure}
\begin{subfigure}[b]{0.49\linewidth}
\centering
\includegraphics[width=0.85\linewidth]{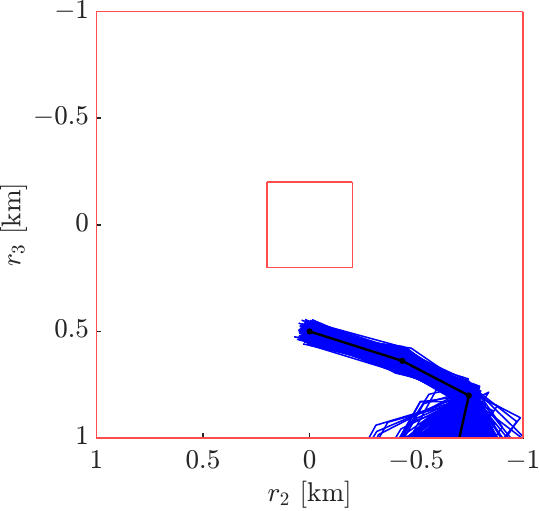}
\end{subfigure}
\caption{Monte Carlo simulation of the NRHO rendezvous maneuver.}
\label{fig:phase1to3-proj-pos}
\end{figure}

\begin{figure}
\centering
\begin{subfigure}[b]{0.49\linewidth}
\centering
\includegraphics[width=\linewidth]{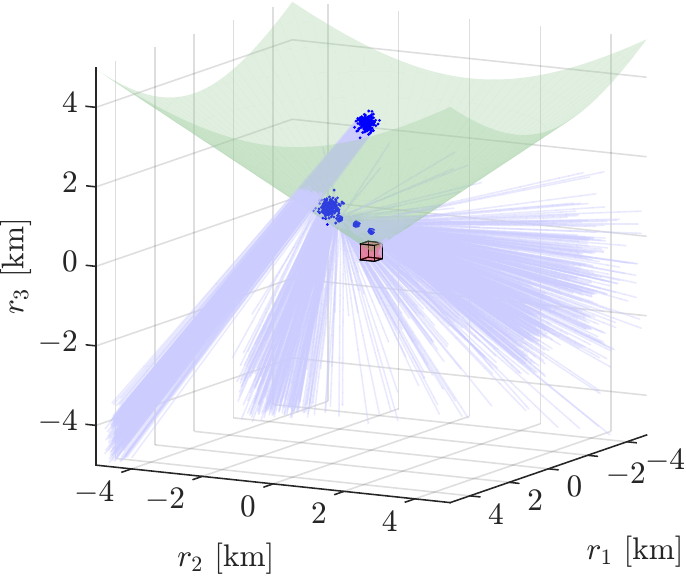}
\caption{Free-drift position trajectories during phase 3 of the rendezvous maneuver.}
\label{fig:phase3-3D-pos}
\end{subfigure}
\hspace{0.1cm}
\begin{subfigure}[b]{0.49\linewidth}
\centering
\includegraphics[width=\linewidth]{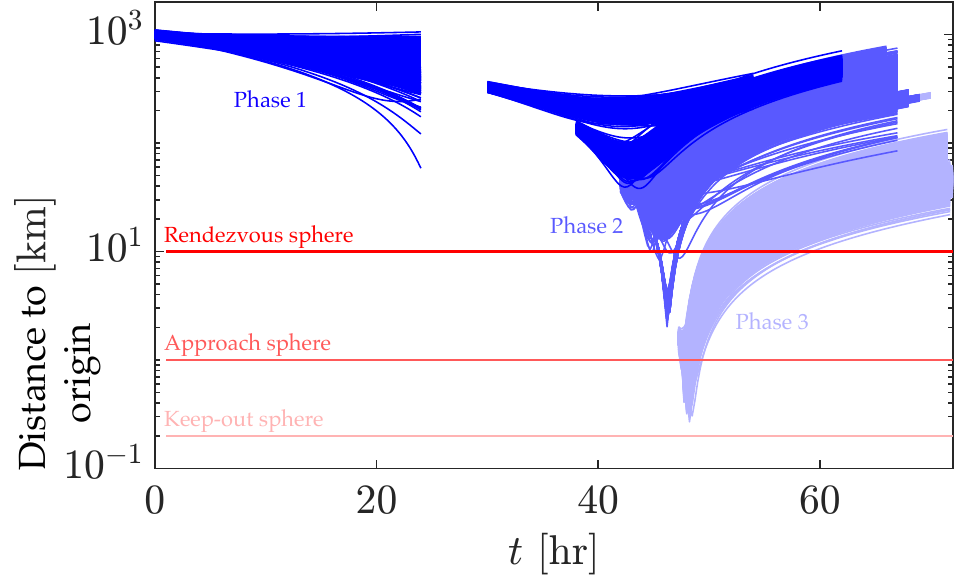}
\caption{Distance of the free-drift position trajectories (blue curves) to the origin.}
\label{fig:dist2origin}
\end{subfigure}
\caption{Free-drift position trajectories of the Monte Carlo samples.}
\label{fig:phase3-3D-pos-dist2origin}
\end{figure}

\begin{figure}[!htp]
\centering
\begin{subfigure}[b]{0.19\linewidth}
\centering
\includegraphics[width=0.85\linewidth]{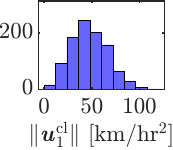}
\end{subfigure}
\begin{subfigure}[b]{0.19\linewidth}
\centering
\includegraphics[width=0.85\linewidth]{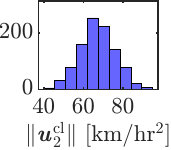}
\end{subfigure}
\begin{subfigure}[b]{0.19\linewidth}
\centering
\includegraphics[width=0.85\linewidth]{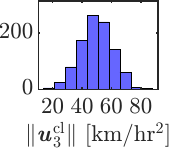}
\end{subfigure}

\vspace{0.5cm}

\begin{subfigure}[b]{0.19\linewidth}
\centering
\includegraphics[width=0.85\linewidth]{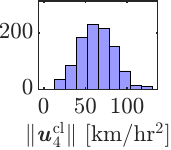}
\end{subfigure}
\begin{subfigure}[b]{0.19\linewidth}
\centering
\includegraphics[width=0.85\linewidth]{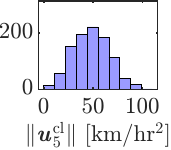}
\end{subfigure}
\begin{subfigure}[b]{0.19\linewidth}
\centering
\includegraphics[width=0.85\linewidth]{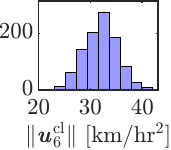}
\end{subfigure}
\begin{subfigure}[b]{0.19\linewidth}
\centering
\includegraphics[width=0.85\linewidth]{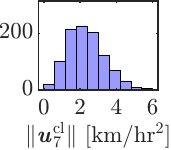}
\end{subfigure}

\vspace{0.5cm}

\begin{subfigure}[b]{0.19\linewidth}
\centering
\includegraphics[width=0.85\linewidth]{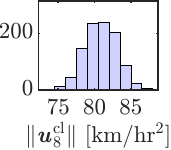}
\end{subfigure}
\begin{subfigure}[b]{0.19\linewidth}
\centering
\includegraphics[width=0.85\linewidth]{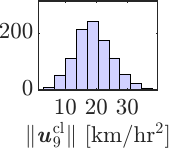}
\end{subfigure}
\begin{subfigure}[b]{0.19\linewidth}
\centering
\includegraphics[width=0.85\linewidth]{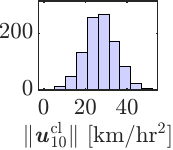}
\end{subfigure}
\begin{subfigure}[b]{0.19\linewidth}
\centering
\includegraphics[width=0.85\linewidth]{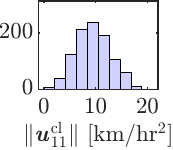}
\end{subfigure}

\vspace{0.5cm}

\begin{subfigure}[b]{0.28\linewidth}
\centering
\includegraphics[width=0.85\linewidth]{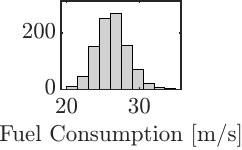}
\end{subfigure}

\vspace{0.5cm}

\begin{subfigure}[b]{\linewidth}
\centering
\includegraphics[width=\linewidth]{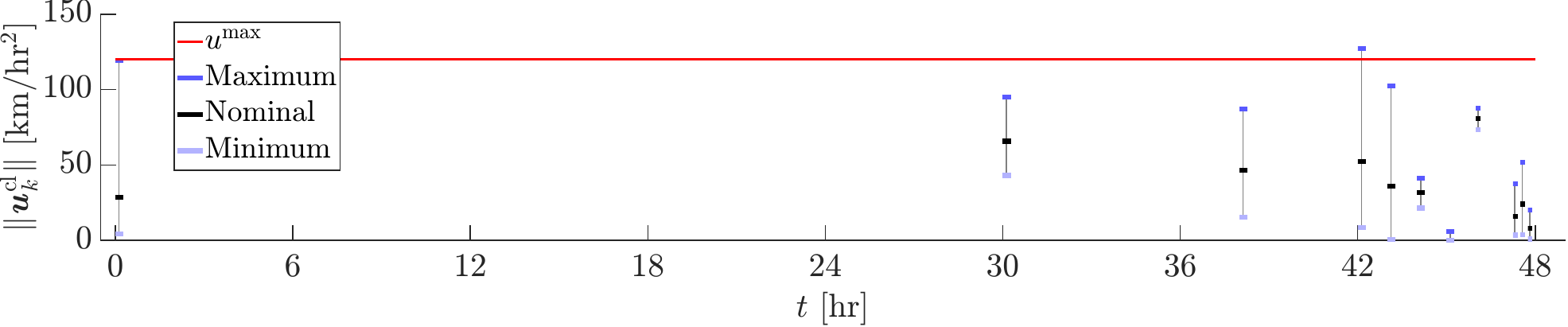}
\end{subfigure}
\caption{Monte Carlo simulation of closed-loop control input magnitudes and total fuel consumption across the three-phase rendezvous maneuver.}
\label{fig:thruster-burn-hist}
\end{figure}

\rtx{The computed rendezvous maneuver and the corresponding Monte Carlo simulation are visualized as follows. Figure~\ref{fig:phase1to3-proj-pos} shows the position trajectories of the Monte Carlo samples in the Sun Referenced LVLH frame (with coordinates $r_1,r_2,r_3$): blue curves denote Monte Carlo samples, black curves denote the nominal trajectory, and dots indicate the time-discretization nodes. The avoid sets for the three phases are shown in red. Figure~\ref{fig:phase3-3D-pos} gives a close-up view of the phase~3 free-drift trajectories at the time-discretization nodes, together with the approach cone (green) and the final avoid set (red). The spacecraft follows a roundabout path to the Gateway to maintain passive safety. Figure~\ref{fig:dist2origin} shows the distance to the origin (the Gateway) for the free-drift position trajectories (blue curves) starting at each thruster firing interval. The darkest blue curve corresponds to phase~1, the medium blue curve to phase~2, and the lightest blue curve to phase~3. The horizontal red lines denote the radii of the rendezvous sphere (10~km; darkest), the approach sphere (1~km; medium), and the keep-out sphere (0.2~km; lightest). The phase~1 distances remain above 10~km, the phase~2 distances remain above 1~km, and the phase~3 distances remain above 0.2~km, indicating that the passive-safety constraint is satisfied in each phase. We do not plot free-drift trajectories along the free-drift arcs between successive thruster firings because they are safe by construction (passive safety is enforced continuously along the trajectory). Figure~\ref{fig:thruster-burn-hist} presents histograms of the closed-loop control input magnitudes across the three-phase maneuver. The first three rows show histograms of control input magnitudes for each phase, and the fourth row shows the histogram of total fuel consumption.}
\begin{table}[!htp]
\centering
\caption{Fuel consumption and constraint satisfaction in 1000-sample Monte Carlo simulation}\label{tab:mc}
\setstretch{1.5}
\begin{tabular}{lc}
\hline\hline
\textbf{Statistic} & \textbf{Value} \\\hline
 Mean~/~Minimum~/~Maximum fuel consumption &  26.3~/~20.2~/~34.8 m/s\\ 
 Samples violating passive-safety constraint & 0 \%\\
 Samples violating approach-cone constraint & 0.8 \%\\
 Samples violating control input constraint & 0.4 \%\\\hline\hline
\end{tabular}
\end{table}

\rtx{A direct comparison between the fuel-consumption statistics from our approach (Table~\ref{tab:mc}) and those from existing NRHO rendezvous methods is challenging because numerical simulations differ in key parameters (Table~\ref{tab:phases-data})}, particularly the spacecraft state at NRHO insertion~\cite{bucchioni2021ephemeris, woffinden2022david, goulet2023robust, nakamura2023rendevous, cavesmith2024angles, cunningham2025robust}. Nonetheless, our method achieves fuel consumption within a comparable range while ensuring continuous-time satisfaction of passive-safety and approach-cone path chance constraints---capabilities not addressed by prior approaches. Additionally, our work is among the few to address \rtx{long-range rendezvous scenarios initiated 1000~km from the Gateway} at NRHO insertion~\cite{nakamura2023rendevous}, while also accounting for underburn safety~\cite{cavesmith2024angles}.

Although the formulation in Section~\ref{sec:prb-form} and the numerical case study address an apolune rendezvous scenario based on published Gateway rendezvous specifications~\cite{nakamura2023rendevous}, where the relative-motion \rtx{dynamics are nearly linear} and the between-maneuver arcs are approximately rectilinear, the approach is not specific to Gateway, apolune, or NRHOs. Beyond requiring well-defined derivatives of the dynamics with respect to state and control input, it is model-based and can be adapted to other rendezvous regimes by updating the application-specific dynamics and constraints, including NRHO perilune rendezvous, \rtx{where nonlinearities are stronger and time scales are shorter~\cite{kikuchi2024extended}}, as well as LEO, GEO, highly elliptical Earth orbits, and orbits about Mars. The core continuous-time chance-constraint framework remains unchanged.
\FloatBarrier
%
\section{Conclusions}\label{sec:conclusions}
We proposed a method for safe, fuel-efficient spacecraft rendezvous with the Gateway using sequential convex programming (SCP). The approach i) enforces passive safety at all times---including underburn cases---and satisfies the approach-cone path constraint in continuous time; ii) meets decision-point requirements along the rendezvous trajectory in accordance with IRSIS guidelines; and iii) models uncertainties from NRHO insertion, actuation errors, and navigation measurements via chance constraints while employing disturbance compensation feedback through correction maneuvers to limit dispersion. To prevent inter-sample violations, we reformulated the continuous-time path constraints as integral constraints using an isoperimetric construction.

\rtx{The formulation improves the efficiency of NRHO rendezvous design and analysis by enabling continuous-time constraint satisfaction without requiring an overly dense temporal discretization for state constraints; only the time grid governing thrust events needs refinement to optimize firing times.} We demonstrated the method on a representative Gateway rendezvous case study and validated it through Monte Carlo simulation. The results show a mean fuel consumption of 26.3~m/s, and each individual constraint is satisfied by more than 99\% of samples---well above the 80\% probability level specified for the chance constraints. Limitations of this work include the omission of certain uncertainty sources (e.g., estimator and actuator biases, and model mismatch) and the absence of a realistic model for the dependence of navigation accuracy on rendezvous trajectory geometry (e.g., relative range, sensor boresight orientation) and sensing modality.

Beyond requiring a dynamics model with well-defined derivatives with respect to state and \rtx{control input}, the approach is dynamics-agnostic: higher-order force models can be incorporated, the fidelity of uncertainty modeling (especially actuation and navigation errors) can be increased, and the disturbance compensation feedback law can be generalized---all without altering the core SCP and continuous-time chance-constraint framework. The method is therefore well suited to safety-critical trajectory optimization problems in which constraints must hold in continuous time under uncertainty, including domains such as autonomous driving~\cite{dicairano2023decision} (with parametrized set-based reachability constraints) and legged locomotion~\cite{liao2023walking} (with continuous-time set-based safety constraints often treated with control barrier functions).%
%
\appendix
\addcontentsline{toc}{section}{Appendix}
\section*{Appendix}
\section{Spacecraft Equations of Motion}\label{app:sc-eom}
The spacecraft equations of motion~\eqref{cislunar-sc-dyn} are expressed in a Moon-centered inertial frame whose axes are aligned with the Earth mean equator and mean equinox of J2000 (EME2000/MEME)~\cite[Sec.~3.7]{vallado2013fundamentals}. In SPICE~\cite{acton1996ancillary}, this corresponds to the \texttt{J2000} orientation with the origin at the Moon's center. The force model includes the two-body Moon term, the lunar $J_2$ perturbation $\aJtwo$, third-body perturbations from the Earth and Sun, and solar radiation pressure (SRP) $\aSRP$. Let $\rsc\in\bR^3$ and $\vsc\in\bR^3$ denote the spacecraft position and velocity in the Moon-centered inertial frame. The nonlinear ODE describing the spacecraft motion is
\begin{align}
    \drsc(t) &= \vsc(t),\label{cislunar-sc-dyn-expand}\\
    \dvsc(t) &= -\,\mu^{\mr{M}} \frac{\rsc(t)}{\|\rsc(t)\|^{3}} \;+\; \aJtwo \;+\; \aSRP \nonumber\\
    &\quad + \mu^{\mr{E}}\!\left( \frac{\rE(t) - \rsc(t)}{\| \rE(t) - \rsc(t) \|^{3}} - \frac{\rE(t)}{\|\rE(t)\|^{3}}\right)
            + \mu^{\mr{S}}\!\left( \frac{\rS(t) - \rsc(t)}{\| \rS(t) - \rsc(t)\|^{3}} - \frac{\rS(t)}{\|\rS(t)\|^{3}}\right),\nonumber
\end{align}
where $z \triangleq (\rsc,\vsc)$ and the right-hand side of~\eqref{cislunar-sc-dyn-expand} defines $\check{f}$ in~\eqref{cislunar-sc-dyn}. The perturbation accelerations $\aJtwo$ and $\aSRP$ are
\begin{subequations}
\begin{align}
    \aJtwo ={} & \rotmatPA(t) \left( -\frac{3\mu^{\mr{M}} J_2 (\Rmoon)^2}{2\|\rsc(t)\|^5} \begin{bmatrix}
        \left(1 - 5\frac{z^{\mr{\scriptscriptstyle PA}}(t)^2}{\|\rsc(t)\|^2}\right)x^{\mr{\scriptscriptstyle PA}}\\[0.15cm]
        \left(1 - 5\frac{z^{\mr{\scriptscriptstyle PA}}(t)^2}{\|\rsc(t)\|^2}\right)y^{\mr{\scriptscriptstyle PA}}\\[0.15cm]
        \left(3 - 5\frac{z^{\mr{\scriptscriptstyle PA}}(t)^2}{\|\rsc(t)\|^2}\right)z^{\mr{\scriptscriptstyle PA}}
    \end{bmatrix} \right)  \\
    \aSRP ={} & P^{\scriptscriptstyle\mr{SRP}}C^{\mr{R}}\Asc (r^{\scriptscriptstyle\mr{AU}})^2 \left(\frac{\rsc(t)-\rS(t)}{\|\rsc(t)-\rS(t)\|^3}\right)
\end{align}\label{J2-SRP}%
\end{subequations}
Among the higher-order lunar gravitational perturbations, only the $J_{2}$ term is retained here \rtx{because} it is the dominant higher-order perturbation~\cite[Fig. 1]{muralidharan2020control}. The $C_{22}$ contribution is also significant but roughly an order of magnitude smaller~\cite{desaedeleer2006combined,li2015autonomous}, and its effect becomes negligible near NRHO apolune. Higher-fidelity models with additional perturbations can be \rtx{incorporated straightforwardly} in the proposed framework. The parameters appearing in~\eqref{cislunar-sc-dyn-expand}--\eqref{J2-SRP} are listed in Table~\ref{tab:cislunar-sc-quantities}. We use the DE421 ephemeris~\cite{folkner2009planetary} within SPICE to query Earth and Sun states relative to the Moon. The SRP model is based on the cannonball model~\cite[Chap. 8]{vallado2013fundamentals}; thus, the illuminated cross section is orientation independent.
\begin{table}[!ht]
\centering
\setstretch{1.5}
\caption{Parameters in the spacecraft equations of motion}
\label{tab:cislunar-sc-quantities}
\begin{tabular}{ll}
\hline\hline
\textbf{Notation} & \textbf{Description} \\\hline
 $\rsc$ & Spacecraft position w.r.t.\ Moon (in Moon-centered J2000 frame)\\
 $\vsc$ & Spacecraft velocity w.r.t.\ Moon (in Moon-centered J2000 frame)\\
 $\rE$  & Earth position w.r.t.\ Moon (in Moon-centered J2000 frame)\\
 $\rS$  & Sun position w.r.t.\ Moon (in Moon-centered J2000 frame)\\
 $\mu^{\mr{M}}$ & Moon gravitational parameter\\
 $\mu^{\mr{E}}$ & Earth gravitational parameter\\
 $\mu^{\mr{S}}$ & Sun gravitational parameter\\
 $\rotmatPA$ & Rotation from Moon principal-axes frame to J2000 frame\\
 $J_2$ & Lunar $J_2$ zonal coefficient ($2.024\times 10^{-4}$)\\
 $\Rmoon$ & Lunar equatorial radius ($1737.1$ km)\\
 $(x^{\mr{\scriptscriptstyle PA}},y^{\mr{\scriptscriptstyle PA}},z^{\mr{\scriptscriptstyle PA}})$ & Spacecraft position in principal-axes frame\\
 $P^{\scriptscriptstyle\mr{SRP}}$ & SRP at 1~AU ($4.5\times 10^{-6}$ Pa)\\
 $r^{\scriptscriptstyle\mr{AU}}$  & Astronomical unit (1~AU)\\
 $C^{\mr{R}}$ & Radiation pressure coefficient (1.15)\\
 $\Asc$ & Area-to-mass ratio of the spacecraft ($16/15000$ m$^2$/kg)\\
\hline\hline
\end{tabular}
\end{table}
\section{Signed Distance}\label{app:signed-distance}
The signed distance~\cite[Chap. IV Sec. 1.3]{hiriart-urruty1993convex} of $z\in\bR^n$ with respect to a nonempty convex set $\mc{D}\subset\bR^n$, denoted by $\rtx{\sd{z}{\mc{D}}}$, is given by
\begin{align}
    \sd{z}{\mc{D}} ={} & \underset{y\in\mc{D}}{\inf}\|z-y\| - \underset{x\notin\mc{D}}{\inf}\|z-x\|.    
\end{align}
The gradient of $\rtx{\sd{\cdot}{\mc{D}}}$ evaluated at $z$, denoted by $\rtx{\nabla_1 \sd{z}{\mc{D}}}$, is given by
\begin{align}
    \rtx{\nabla_1 \sd{z}{\mc{D}}} ={} & \frac{(z-\projboundary{z}{\mc{D}})^\top\!\!\!\!}{\rtx{\sd{z}{\mc{D}}}},    
\end{align}
where $\projboundary{z}{\mc{D}} = \underset{y\in\partial\mc{D}}{\mr{argmin}}\|z-y\|$ is the projection of $z$ onto the boundary of $\mc{D}$, denoted by $\partial\mc{D}$. 

Consider a polyhedron of the form $\mc{D} = \{z\in\bR^n\,|\,Hz\le h\}$, where $H = [H_1\,\ldots\,H_m]^\top$ and $h=(h_1,\ldots,h_m)$, with $H_i\in\bR^n$ and $h_i\in\bR$ for $i=1,\ldots,m$. Then,
\begingroup
\everymath{\displaystyle}
\begin{align}
    \rtx{\sd{z}{\mc{D}}} ={} \left\{ \begin{array}{ll} \underset{y\in\mc{D}}{\min}\|z-y\| &~\text{if}~z\notin\mc{D}\\-\underset{1\le i \le m}\min\frac{|H_i^\top z-h_i|}{\|H_i\|}&~\text{otherwise}\end{array} \right.
\end{align}
\endgroup
where the first case amounts to solving an inequality-constrained QP and the second case is a simple algebraic enumeration \rtx{over at most all} faces of the polyhedron. Furthermore, the projection of $z$ onto the boundary of $\mc{D}$ is given by
\begingroup
\everymath{\displaystyle}
\begin{align}
    \projboundary{z}{\mc{D}} ={} \left\{ \begin{array}{ll} \underset{y\in\mc{D}}{\mr{argmin}}\|z-y\| &~\text{if}~z\notin\mc{D}\\z - \frac{(H_{i^\star}^\top z - h_{i^\star})}{\|H_{i^\star}\|^2}H_{i^\star} &~\text{otherwise}\end{array} \right.
\end{align}
\endgroup
where ${i^\star} = \underset{1\le i \le m}{\mr{argmin}}|H_i^\top z - h_i|/\|H_i\|$.
%
%
\bibliography{references-revised}
\end{document}